\documentclass[12pt]{amsart}

\usepackage{amsmath}
\usepackage{amssymb}
\usepackage{enumerate}
\usepackage{euscript}
\usepackage{amscd}
\usepackage[all]{xy}

\newtheorem{thm}{Theorem}[section]
\newtheorem{lem}[thm]{Lemma}
\newtheorem{prop}[thm]{Proposition}
\newtheorem{cor}[thm]{Corollary}
\newtheorem{defn}[thm]{Definition}

\newtheorem{clm}{Claim}

\numberwithin{equation}{section}

\newcommand{\A}{\mathbb A}
\newcommand{\Z}{\mathbb Z}
\newcommand{\F}{\mathbb F}
\newcommand{\C}{\mathbb C}
\newcommand{\Q}{\mathbb Q}
\newcommand{\R}{\mathbb R}

\newcommand{\B}{\mathcal B}
\newcommand{\K}{\mathcal K}
\newcommand{\T}{\mathcal T}
\newcommand{\cd}{\mathcal D}
\newcommand{\ck}{\mathcal K}
\newcommand{\co}{\mathcal O}
\newcommand{\cs}{\mathcal S}
\newcommand{\Sc}{\mathcal S}

\newcommand{\sB}{\EuScript{B}}
\newcommand{\Cs}{\EuScript{C}}
\newcommand{\Ss}{\EuScript{S}}
\newcommand{\E}{\EuScript{E}}

\newcommand{\gD}{\mathfrak D}

\newcommand{\gM}{{\mathfrak M}}
\newcommand{\gC}{\mathfrak C}
\newcommand{\gc}{\mathfrak c}
\newcommand{\gp}{\mathfrak p}

\newcommand{\esN}{{\EuScript{N}}}
\newcommand{\esR}{{\EuScript{R}}}

\newcommand{\rew}{{\rm Re}(w)}
\newcommand{\rez}{{\rm Re}(z)}
\newcommand{\res}{{\rm Re}(s)}

\def\ord{\mathop{\mathrm{ord}}}
\def\vol{\mathop{\mathrm{vol}}}
\def\Res{\mathop{\mathrm{Res}}}
\def\supp{\mathop{\mathrm{supp}}}

\newcommand{\adl}{\mathbb A}
\newcommand{\idl}{\mathbb A^\times}
\newcommand{\adlint}[1]{\int_{\adl/k}\hspace{-3mm}d{#1}}
\newcommand{\idllength}{-6mm}
\newcommand{\idlint}[1]{\int_{\idl/k^\times}\hspace{\idllength}d^\times{#1}}

\newcommand{\idn}[1]{|{#1}|_\adl}

\newcommand{\map}{\rightarrow}

\newcommand{\genus}{\mathfrak g}
\newcommand{\vprod}{\prod_{v\in \gM}}
\newcommand{\cmpint}{\frac{1}{2\pi i}\int_{\rez=r}}
\newcommand{\ac}[1]{\langle{#1}\rangle}

\newcommand{\nd}[1]{\esN(\gD_{{#1}})}
\newcommand{\sr}{{\rm SR}}
\newcommand{\vep}{\varepsilon}

\newcommand{\stwtw}[4]
{\left(\begin{smallmatrix}{#1}&{#2}\\{#3}&{#4}\\\end{smallmatrix}\right)}
\newcommand{\twtw}[4]
{\left(\begin{array}{cc}{#1}&{#2}\\{#3}&{#4}\\\end{array}\right)}
\newcommand{\svc}[2]
{\left(\begin{smallmatrix}{#1}\\{#2}\\\end{smallmatrix}\right)}

\makeatletter

\renewcommand{\subsection}{\@startsection{subsection}{2}%
  \z@{.5\linespacing\@plus.7\linespacing}{-.5em}%
  {\normalfont\bfseries\S\,}}

\makeatother

\begin{document}

\title[mean value theorem]
{A mean value theorem for orders of
degree zero divisor class groups of
quadratic extensions\\
over a function field}

\author[Takashi Taniguchi]{Takashi Taniguchi}
\address{Department of Mathematical Sciences, University of Tokyo}
\email{tani@ms.u-tokyo.ac.jp}
\date{\today}
\begin{abstract}
Let $k$ be a function field of one variable over a finite field
with the characteristic not equal to two.
In this paper, we consider the prehomogeneous representation
of the space of binary quadratic forms over $k$.
We have two main results.
The first result is on the principal part of the global zeta function
associated with the prehomogeneous vector space.
The second result is on a mean value theorem
for degree zero divisor class groups of
quadratic extensions over $k$,
which is a consequence of the first one.

\end{abstract}

\maketitle
\setcounter{tocdepth}{1}

\section{Introduction}

Let $k$ be a function field of one variable over a finite field
and $h_k$ be the order of degree zero divisor class group.
We assume that ${\rm char}(k)\not=2$.
For a quadratic extension $L$ of $k$,
we denote by $\nd{L}$
the norm of its relative discriminant.
In this paper, we will give a mean value of
$h_L$ with respect to $\nd{L}$.

Our main result is Theorem \ref{thm:mv}.
We briefly state our results here.
Let $q$ and $\zeta_k(s)$ be the order of the constant field and
the Dedekind zeta function of $k$.
We denote by $\gM$ the complete set of places of $k$,
and by $q_v$ the absolute norm of $v\in\gM$.
The following theorem is a special case of Theorem \ref{thm:mv}.
\begin{thm}\label{thm:intro}
\begin{equation*}
\lim_{n\to\infty}\frac{1}{q^{3n}}
\sum_{\substack{[L:k]=2\\ \nd{L}=q^{2n}}}h_{L}
=2h_k\frac{h_k}{q-1}\zeta_k(2)\prod_{v\in\gM}(1-q_v^{-2}-q_v^{-3}+q_v^{-4}).
\end{equation*}
\end{thm}

Theorems of this kind are called density theorems.
Today, many density theorems are known.
Among them,
theorems about the asymptotic behavior of the mean value
of the number of equivalence classes of
quadratic forms are very classical
and studied by many mathematicians including Gauss, Siegel, and Shintani.
The density theorem for the case of binary quadratic forms is known as
Gauss' conjecture. 
This was firstly proved
by Lipschitz for imaginary case,
and by Siegel for real case.
Siegel \cite{si} also proved
the density theorem for integral equivalence classes
of quadratic forms in general.

M.Sato and T.Shintani formulated this kind of density problems
using the notion of prehomogeneous vector spaces.
In \cite{Shi2}, Shintani considered
a zeta function associated with
the space of quadratic forms.
There, he reproved the Gauss' conjecture, 
and improved the error estimate.
Shintani \cite{Shi1} also considered the space of binary cubic forms,
and gave the density of the class number of integral binary cubic forms.

Datskovsky and Wright
\cite{Wr}, \cite{Dat1}, \cite{Da-Wr1}, \cite{Da-Wr2}
treated the zeta function associated with
the space of binary cubic forms over an arbitrary global field
using the adelic language.
Then they obtained the Davenport-Heilbronn \cite{dahe}, \cite{dahe2}
density theorem of cubic fields
from the viewpoint of prehomogeneous vector spaces.

Also for the space of binary quadratic forms,
Shintani's global theory was extended to
an arbitrary number field by Yukie \cite{yu1}.
Then Datskovsky \cite{Dat1} carried out the local theory,
and gave an another proof of the 
Goldfeld-Hoffstein \cite{goho} mean value theorem of
``class number times regulator'' of quadratic extensions
over any fixed number field.
However, the problem for
binary quadratic forms over a function field
was left unanswered. In this paper, we will study this case.

For the general process
from a prehomogeneous vector space
to its density theorem,
see \cite{yu2} or \cite{yu}.

This paper is organized as follows.
In Section \ref{sec:nt},
we collect notation we use throughout this paper.
In Section \ref{sec:22},
we review structures of the space of binary quadratic forms.
The rational orbit decomposition and
the structure of the stabilizers are discussed.
In Sections \ref{sec:g1} and \ref{sec:g2}, we treat the global theory.
In Section \ref{sec:g1},
we define the global zeta function $Z(\Phi,s)$.
This is a function of
a Schwartz-Bruhat function $\Phi$,
and a complex variable $s$.
The main purpose of this section is to show
that it converges if $\res$ is sufficiently large.
In Section \ref{sec:g2},
we study analytic properties of the global zeta function.
We show its rationality and determine the pole structure.
Our main result is Theorem \ref{th:gz}.
The residues at the poles are described by means of some distributions.
We later need its rightmost pole, which is rather simple;
it is a constant multiple of $\hat\Phi(0)$
where $\hat\Phi$ is the Fourier transform of $\Phi$.

In Section \ref{sec:l}, we define various invariant measures,
and consider local zeta functions.
Then in Section \ref{sec:mv}, 
we define certain Dirichlet series and study their analytic properties
by putting together the results we have obtained before.
And by using the filtering process, we obtain a density theorem
which is a generalized version of Theorem \ref{thm:intro}.

Our discussion is quite similar to \cite{Dat1},
which is the case of number fields.
However, as A.Yukie pointed out at \cite{yu2}, p.12, there is an
incomplete argument in his paper.
His choice of the measure on the stabilizer at \cite{Dat1}, p.218 is wrong
because it does not satisfy the functorial property,
which he implicitly used in \cite{Dat1}, p.230.
We will correct the argument in Section \ref{sec:l} following
\cite{kayu}.
It is easy to see that
our choice can also be applied to his paper,
and the final results of \cite{Dat1} need no modification.

\vspace{0.2in}

\noindent
{\bf Acknowledgments.} 
The author would like to express
his sincere gratitude
to his advisor T.Terasoma
for the support and encouragement.
The author is also deeply grateful to A.Yukie,
who read the manuscript
and gave many helpful suggestions.

\section{Notation}\label{sec:nt}

Here, we will prepare basic notation.
For a finite set $X$ we denote by $\#X$ its cardinality.
If $f,g$ are functions on a set $Z$, and $|f(z)|\leq Cg(z)$
for some constant $C$ independent of $z\in Z$,
we denote $f(z)\ll g(z)$.
The standard symbols $\Q,\R,\C,$ and $\Z$ will denote the set of
rational, real, complex numbers and the rational integers, respectively.
If $R$ is any ring then $R^\times$ is the group of units of $R$,
and if $X$ is a variety defined over $R$ then $X_R$ denotes its $R$-points.

Suppose that $G$ is a locally compact group
and $\Gamma$ is a discrete subgroup of $G$
contained in the maximal unimodular subgroup of $G$.
For any left invariant measure $dg$ on $G$,
we choose a left invariant measure $dg$
(we use the same notation,
but the meaning will be clear from the context)
on $X=G/\Gamma$ so that
\begin{equation*}
\int_G f(g)dg=\int_X\sum_{\gamma\in\Gamma}f(g\gamma)dg.
\end{equation*}

Throughout this paper, we agree that
$k$ denotes any fixed function field of one variable
over a finite field of constants $\F_q$,
$q\not =2^n$.
Denote by $\gM$ the complete set of places of $k$.
For $v\in\gM$, $k_v$ denotes the completion of $k$ at $v$ and
$|\cdot|_v$ the normalized absolute value on $k_v$.
We denote by $\co_v, \pi_v$ and $q_v$,
the ring of integers of $k_v$,
a fixed uniformizer in $\co_v$,
and the cardinality of $\co_v/\pi_v\co_v$, respectively.

Returning to $k$, let $\genus,h_k$ be the genus of $k$,
the order of degree zero divisor class group of $k$, respectively.
It will be convenient to set
$$
\gC_k=\frac{h_k}{q-1}.
$$
Let $\zeta_k(s)$ denote the Dedekind zeta function 
$\prod_{v\in\gM}(1-q_v^{-1})^{-s}$ of $k$.
It is known that $\zeta_k(s)$ is a rational function of $q^s$ and
$$\lim_{s\rightarrow 1}(1-q^{1-s})\zeta_k(s)=q^{1-\genus}\gC_k.$$

Later, we will consider quadratic extensions of $k$.
For an extension $L$ of $k$,
the symbols $q_L, h_L, \gC_L$ and $\gM_L$
are defined in similar way.
We will denote by $\nd{L}$ 
the ideal norm of the relative discriminant of $L$ over $k$.

Let $\A,\A^\times$ be
the ring of adeles and the group of ideles of $k$, respectively.
The group $\A^\times$ is endowed with the idele norm,
denoted by $\idn{\cdot}$, 
where $\idn{x}=\vprod |x_v|_v,x=(x_v)_v\in\A^\times$.
The field $k$, identified with a subset of $\A$ via the diagonal embedding,
forms a lattice in $\A$.
Let $\A^\times\supset \A^n$ be the set of elements with idele norm $q^n$.
If $V$ is a vector space over $k$,
let $\cs(V_\A),\cs(V_{k_v})$
be the space of Schwartz-Bruhat functions.

We choose a Haar measure
	$dx$			on	$\adl$,
	$d^\times t$	on	$\adl^\times$,
	$dx_v$			on	$\co_v$,
	$d^\times t_v$	on	$\co_v^\times$,
so that
$$
	\int_{\adl/k} dx = 1,					\quad
	\int_{\adl^0/k^\times} d^\times t=1,	\quad
	\int_{\co_v}dx_v=1,						\quad
	\int_{\co_v^\times} d^\times t_v=1,
$$
respectively.
We later have to 
compare the global measure and the product of 
local measures.
It is well known that
\begin{equation}\label{productmeasures} 
dx = q^{1-\genus}\prod_{v\in\gM} dx_v\qquad \text{and}\qquad 
dt = \gC_k^{-1}\prod_{v\in\gM} d^{\times} t_v.
\end{equation}

Finally, for the sake of convenience ,
we denote the element $G=GL(1)\times GL(2)$ in the following manner:
\begin{equation*}
n(u)=		\left(1,\twtw 10u1\right),\ \ \ \ 
d(t,t_1)=	\left(t,\twtw{t_1}00{t_1}\right),\ \ \ \ 
a(\tau)=	\left(1,\twtw 100\tau\right).
\end{equation*}
We note that $a(\tau)n(u)=n(\tau u)a(\tau)$.

\section{The space of binary quadratic forms}\label{sec:22}

Let $V$ be the 3-dimensional affine space.
We define $V$ with the space of binary quadratic forms
via the correspondence:
$$
x=(x_0,x_1,x_2)\in V\longleftrightarrow
F_x(z_1,z_2)=x_0z_1^2+x_1z_1z_2+x_2z_2^2.
$$

The group $G=GL(1)\times GL(2)$ acts on $V$ from the left,
the $GL(1)$-part acts on $V$ by the usual scalar multiplication,
and the $GL(2)$-part acts on $V$ by the linear change of variables.
Explicitly, the action of $G$ on $V$ is given by
$$
F_{g\cdot x}(z_1,z_2)=tF_x(az_1+cz_2,bz_1+dz_2),\quad
g=\left(t,\twtw abcd\right)\in G, x\in V.
$$

For any field $K$ the space $V_K$ is a 3-dimensional vector space over $K$.
The action of $G$ on $V$ then gives a representation
$\varrho:G\rightarrow GL(V)$ defined over $K$.
The kernel of $\varrho$ is a one-dimensional torus 
$T_\varrho$ in the center of $G$.
Let
$$H={\rm Im}(\varrho).$$
$H$ is a closed reductive subgroup of $GL(V)$.

For $x\in V$, let $P(x)$ denote the discriminant of $F_x$:
$$
P(x)={x_1}^2-4x_0x_2.
$$
For $g=(t,\stwtw abcd)\in G$, set $\chi(g)=t(ad-bc)$.
Since $\chi(T_\varrho)=1$, we can regard $\chi$ as a character on $H$.
It is easy to see that
$$P(h\cdot x)=(\chi(h))^2P(x) \quad \text{and} \quad
\det(h)=\chi(h)^3.$$

We call a form $x$ non-singular if $P(x)\not =0$
and singular otherwise.
Let $V'$ be the set of all non-singular forms in $V$.
It is easy to see that
two forms in $V_K'$ are $G_K$-equivalent
if and only if their splitting fields over $K$ are the same.
Thus non-singular $G_K$-orbits in $V_K$
are in one-to-one correspondence with the extensions of $K$
of degree less or equal to two.
For $x\in V_K$, $K(x)$ denotes the splitting field of $x$ over $K$.
And for $x\in V_K'$, let $G_x\subset G$ denote the stabilizer subgroup of $x$,
and $G_x^\circ$ the identity component.

\begin{prop}\label{prop:Gx0}
\begin{enumerate}[{\rm (i)}]
\item $|G_x/G_x^\circ|=2$.
\item
$G_x^\circ=
	\begin{cases} 
		GL(1) \times GL(1) & K(x)=K,\\
		G_x^\circ=R_{K(x)/K}(GL(1)) & [K(x):K]=2.
	\end{cases}$
\end{enumerate}
\end{prop}

\proof

Here, we write the summary which will be needed in Section \ref{sec:l}.
For detail, see \cite{Dat1}.

First for the case $K(x)=K$.
Let $F_x(z_1,z_2)=(pz_1+qz_2)(rz_1+sz_2)$,
then for $g=(t,\stwtw abcd)\in G_x^\circ$,
$\stwtw abcd$ acts on $pz_1+qz_2,rz_1+sz_2$ as scalar multiplications,
hence it has two eigenvectors and they are $\svc pq,\svc rs$.
Set these eigenvalues $\alpha,\beta$,
then $t=(\alpha\beta)^{-1}$ and
$$
G_x^\circ \rightarrow GL(1)\times GL(1)
\qquad g\mapsto (\alpha,\beta)
$$
gives an isomorphism of groups.
Note that $g\mapsto (\beta,\alpha)$ also gives an isomorphism.

Next for the case $[K(x):K]=2$.
Let
$$
F_x(z_1,z_2)=x_0z_1^2+x_1z_1z_2+x_2z_2^2
=x_0(z_1+\theta z_2)(z_1+\theta' z_2),
$$
where $\theta'$ is the Galois conjugate of $\theta$ over $K$.
Then for $g=(t,\stwtw abcd)\in G_x^\circ$,
$\stwtw abcd$ acts on $z_1+\theta z_2$ as scalar multiplication by $a+b\theta$
and $t=N_{K(x)/K}(a+b\theta)^{-1}$.
Now, let
$$
\varphi:G_x^\circ \rightarrow GL(1)
\qquad g\mapsto a+b\theta.
$$
Then $\varphi$ is a morphism of linear algebraic groups,
defined over $K(x)$, and as an abstract group,
$G_x^\circ(K)\cong GL(1)_{K(x)}$.
Also note that $g\mapsto a+b\theta'$ gives another isomorphism.
\hspace*{\fill}$\square$

\section{The global zeta function: definition and convergence}\label{sec:g1}
\subsection{Definition of the global zeta function}

From now on, 
let $k$ be a function field of one variable
over a finite field of constants $\F_q$ such that ${\rm char}(k)\not=2$.
Then $H_k$ becomes a discrete subgroup of $H_\A$
and $V_k$ a $H_k$-invariant lattice of $V_\A$.

\begin{defn}
Let $V_k''=\{x\in V_k'\mid[k(x):k]=2\}$.
For $\Phi\in\Sc(V_\A)$ and $s\in\C$,
we define the global zeta function $Z(\Phi,s)$ by
$$
Z(\Phi,s)=\int_{H_\A/H_k}\idn{\chi(h)}^s
\sum_{x\in V_k''}\Phi(h\cdot x)dh
$$
where $dh$ is a Haar measure on $H_\A$
that will be normalized in \S\ref{subsec:dh1}.
\end{defn}

\subsection{A Haar measure on $H_\adl$}\label{subsec:dh1}

Firstly, we describe a Haar measure on $G_\A$.
Let $G_\adl=\K B_\adl$ be the Iwasawa decomposition,
where $\K$ is the standard maximal compact subgroup of $G_\adl$
and $B_\A$ a Borel subgroup of $G_\adl$.
More precisely,
$\K=\prod_{v\in \gM}\K_v, \K_v=G_{\co_v}$, and
$B=\{(t,\stwtw {t_1}0a{t_2})\}$.

Every element of $B_\adl$ can be written uniquely as
$b=d(t,t_1)a(\tau)n(u)$, where $u\in\adl,t,t_1,\tau\in\idl$,
and it is easy to check that 
$db=\idn{\tau}d^\times td^\times t_1 d^\times\tau du $
is a right invariant measure on $B_\adl$. 
We normalize the Haar measure $d\kappa$ on $\K$
by $\int_\K d\kappa=1$.
Then $dg=d\kappa db$ gives a normalization
of the Haar measure on $G_\adl$.

Let $\A(\emptyset)=\prod_{v\in \gM}\co_v,
\A^\times(\emptyset)=\prod_{v\in \gM}\co_v^\times$.
Note that $t,t_1,u,\tau$ of the Iwasawa decomposition
$g=\kappa d(t,t_1)a(\tau)n(u)$ are not unique,
but $t,t_1,\tau$ are unique up to multiplication
by elements of $\idl(\emptyset)$
(hence $\idn{t(g)},\idn{t_1(g)},$ and $\idn{\tau(g)}$
are well-defined),
and $u$ is uniquely determined modulo $\tau^{-1}\A(\emptyset)$.

Recall that $H\cong G/T_\varrho$
where $T_\varrho=\{d(t_1^2,t_1^{-1})\in G\}$.
Define $dh$ on $H_\A$ by setting $dg=d^\times t_1dh$.
More explicitly, if we write
$$h=\varrho(\kappa d(t,1)a(\tau)n(u)),$$
then
$$dh=|\tau|_\A d\kappa d^\times t d^\times\tau du.$$

Let $\B=\varrho(B)$ and $\T=\varrho(T)$
where $T$ is the maximal torus in $B$.
In Section \ref{sec:g2}, we will compute integrals
over $H_\adl/\B_k$ and $H_\adl/\T_k$
with respect to the measure $H_\A$.
Such integrals are explicitly given as follows:

\begin{align}
\int_{H_\adl/\B_k}\Psi(h)dh 
&=	\int_\K d\kappa\idlint{t}\idlint{\tau}\adlint{u}\ \ 
	\Psi(\varrho(\kappa d(t,1)a(\tau)n(u)))\idn{\tau},
\label{eq:edh0}\\
\int_{H_\adl/\T_k}\Psi(h)dh 
&=	\int_\K d\kappa\idlint{t}\idlint{\tau}\int_\adl du\ \ 
	\Psi(\varrho(\kappa d(t,1)a(\tau)n(u)))\idn{\tau},
\label{eq:edh1}
\end{align}
where $\Psi\in L^1(H_\adl/\B_k),L^1(H_\adl/\T_k)$, respectively.

Let $\widehat{\T}_k$ be the subgroup of $H_k$
generated by $\T_k$ and $\iota=\varrho(1,\stwtw 0110)$.
We will also consider an integral over $H_\A/\widehat{\T}_k$.
Note that $\T_k$ is a subgroup
of $\widehat{\T}_k$ of index 2.
In the same way as \cite{yu1}, p362, we can prove the next lemma.
\begin{lem}
For $u=(u_v)_{v\in \gM}\in\A$, define
$$
\alpha(u)=\prod_{v\in \gM}\max(1,|u_v|_v).
$$
Then for $g=\kappa d(t,1)n(u)a(\tau)\in G_\adl$,
$\idn{\tau(g\iota)}=\alpha(u)^2\idn{\tau(g)}^{-1}$.
\end{lem}

By this lemma, 
we can write down the integral over $H_\adl/\widehat{\T}_k$
as follows:
\begin{equation}
\int_{H_\adl/\widehat{\T}_k}\Psi(h)dh
=	\int_\K d\kappa\idlint{t}
	\left(
		\underset{\idn{\tau}<\alpha(u)}{\idlint{\tau}}
		+\frac{1}{2}\underset{\idn{\tau}=\alpha(u)}{\idlint{\tau}}
	\right)
	\int_{\adl}du\ \ 
	\Psi(\varrho(\kappa d(t,1)n(u)a(\tau))),
\label{eq:edh2}
\end{equation}
where $\Psi\in L^1(H_\adl/\widehat{\T}_k)$.

\subsection{The convergence of the global zeta function}\label{subsec:cnv}

Firstly, we describe
a fundamental domain of $H_\A/H_k$.
Set $\Theta_2=\{g\in GL(2)_\A|\idn{\tau(g)}<q^{2\genus}\}$ and
$\Theta=\idl\times\Theta_2\subset G_\adl$.
Then by a reduction theorem of Harder \cite{Ha},
$\Theta G_k=G_\A$.
The set $\Theta$ is invariant under multiplication
by elements of $B_k$ on the right.
Hence a fundamental domain of $B_k$ in $\Theta$ contains
a fundamental domain of $G_k$ in $G_\A$.
Then
$$
|Z(\Phi,s)|\ll
\int_{\varrho(\Theta)/\B_k}\idn{\chi(h)}^{\res}
\sum_{x\in V_k''}|\Phi(h\cdot x)|dh.
$$

Let us describe
the fundamental domain of $\varrho(\Theta)/\B_k$.
\cite{Dat2} Lemma 2.2 immediately leads to the following.
\begin{lem}\label{lem:Si}
Every element of $\varrho(\Theta)$ is 
right $\B_k$-equivalent to an element of
 $\varrho(\Ss)$, where
$$
\Ss=\bigcup_{u}\bigcup_{t,\tau}
\K n(u)d(t,1)a(\tau)\subset G_{\adl},
$$
$t,\tau$ run over a set of 
representative of $\A^\times/k^\times$,
$\idn{\tau}\leq q^{2\genus}$,
and $u$ runs over a finite set in $\A$.
\end{lem}

Next, we will estimate 
$\sum_{x\in V_k''}\Phi(h\cdot x)$.

\begin{lem}\label{lem:t}
There exists an integer $N>0$ such that
$\sum_{x\in V_k''}|\Phi(d(t,1)a(\tau)\cdot x)|=0$
whenever $\idn{t}>q^N$.

\end{lem}
\begin{proof}
Let $x=(x_0,x_1,x_2)\in V_k''$. Suppose
$$
\varrho(d(t,1)a(\tau))\cdot x=
\left(t x_0,t\tau x_1,t\tau^2 x_2\right)\in\supp(\Phi).
$$
Since $\supp(\Phi)$ is compact, the first coordinate is bounded,
that is, there exists an integer $N>0$ such that
$\idn{t x_0}\leq q^{N}$.
On the other hand, $x_0\not=0$ follows from $x\in V_k''$,
and hence $\idn{x_0}=1$.
This completes the proof.
\end{proof}

We recall the following well known fact.
\begin{lem}\label{lem:weil}
\begin{enumerate}[{\rm (i)}]
\item
Let $\Psi$ be a Schwartz-Bruhat function on $\A$.
Then,
$$
\sum_{x\in k}|\Psi(tx)|\ll\max(1,|t|_\adl^{-1}),\quad
\sum_{x\in k^\times}|\Psi(tx)|\ll|t|_\adl^{-1}.
$$
\item
Let $\Psi$ be a Schwartz-Bruhat function on ${\rm Aff}^n_\A$.
Then there exist Schwartz-Bruhat functions $\Psi_1,\ldots,\Psi_n\geq0$
such that
$$|\Psi(x_1,\ldots,x_n)|\leq \Psi_1(x_1)\cdots\Psi_n(x_n)$$
for $x_1,\ldots,x_n$.
\end{enumerate}
\end{lem}

We are now ready to prove

\begin{thm}\label{thm:cnv}
$Z(\Phi,s)$ converges absolutely and locally uniformly
in the region of $\res>3$.
In particular, $Z(\Phi,s)$ is a holomorphic function of $s$ in the region.
\end{thm}
\begin{proof}
Let $\sigma=\res$. We have
\begin{align*}
Z(\Phi,s)
&\ll	\idlint{t}\underset{|\tau|_\adl\leq q^{2\genus}}{\idlint{\tau}}\ \ 
		\idn{\tau}|t\tau|_\adl^{\sigma}
		\sum_{x\in V_k''}\left|\Phi(tx_0,t\tau x_1,t\tau^2x_2)\right|
			\qquad \text{by Lemma \ref{lem:Si}}\\
&=		\underset{|t|_\adl\leq q^N}{\idlint{t}}
		\underset{|\tau|_\adl\leq q^{2\genus}}{\idlint{\tau}}\ \ 
		\idn{\tau}|t\tau|_\adl^{\sigma}
		\sum_{x\in V_k''}\left|\Phi(tx_0,t\tau x_1,t\tau^2x_2)\right|
			\qquad \text{by Lemma \ref{lem:t}}.\\
\intertext{For $x=(x_0,x_1,x_2)\in V_k''$,
both $x_0\not=0$ and $x_2\not=0$.
Hence,}
&\ll	\underset{|t|_\adl\leq q^N}{\idlint{t}}
		\underset{|\tau|_\adl\leq q^{2\genus}}{\idlint{\tau}}\ \ 
		\idn{\tau}|t\tau|^{\sigma}|t|^{-1}|t\tau^2|^{-1}
		\max(1,|t\tau|_\adl^{-1})
			\qquad \text{by Lemma \ref{lem:weil}}\\
&=		\underset{|t|_\adl\leq q^{N+2\genus}}{\idlint{t}}
		\underset{|\tau|_\adl\leq q^{2\genus}}{\idlint{\tau}}\ \ 
		\idn{\tau}|t|_\adl^{\sigma}
		\idn{t/\tau}^{-1}\idn{t\tau}^{-1}
		\max(1,\idn{t}^{-1})\\
&\ll	\underset{|t|_\adl\leq q^{N+2\genus}}{\int_{\idl/k^\times}}
		\idn{t}^{\sigma-3}\ d^\times t
		\cdot
		\underset{|\tau|_\adl\leq q^{2\genus}}{\int_{\idl/k^\times}}
		\idn{\tau}\ d^\times \tau.
\end{align*}
Hence, we obtain the desired result.
\end{proof}

\subsection{The definition of $Z_+(\Phi,s)$ and $Z_-(\Phi,s)$}

For $t\in \idl$, set
$$\lambda_+(t)=
	\left\{\begin{array}{ll}
		0, 		& \idn{t}<1,\\
		\frac{1}{2}, 	& \idn{t}=1,\\
		1, 		& \idn{t}>1,\\
	\end{array}\right.\quad
\lambda_-(t)=
	\left\{\begin{array}{ll}
		1, 		& \idn{t}<1,\\
		\frac{1}{2}, 	& \idn{t}=1,\\
		0, 		& \idn{t}>1,\\
	\end{array}\right.$$
and for $h\in H_\adl$, the same symbol $\lambda_+,\lambda_-$ denotes
$$
\lambda_+(h)=\lambda_+(\chi(h))\quad \lambda_-(h)=\lambda_-(\chi(h)).
$$
Using these symbols, define
\begin{align*}
Z_+(\Phi,s)
&=	\int_{H_\adl/H_k}\idn{\chi(h)}^s\lambda_+(h)
	\sum_{x\in V_k''}\Phi(h\cdot x)dh,\\
Z_-(\Phi,s)
&=	\int_{H_\adl/H_k}\idn{\chi(h)}^s\lambda_-(h)
	\sum_{x\in V_k''}\Phi(h\cdot x)dh.
\end{align*}
Then we have
$$Z(\Phi,s)=Z_+(\Phi,s)+Z_-(\Phi,s).$$
Moreover, the following proposition holds.
We can obtain this just as in \cite{Dat1} Proposition 2.1,
and we omit the proof.
\begin{prop}\label{prop:z+}
$Z_+(\Phi,s)$ is a polynomial in $q^s$.
\end{prop}

\section{The global zeta function: analytic continuation 
and the principal part formula}\label{sec:g2}

\subsection{The principal part $I(\Phi,s)$}
Let $\ac{\ }:\adl\map\C^\times$
be a fixed non-trivial additive character on $\A$
trivial on $k$.
Let $[\ ,\ ]$ be a nondegenerate symmetric bilinear form
on $V_\A$ given by
$$
[x,y]=x_0y_2-\frac{1}{2}x_1y_1+x_2y_0,
$$
and we identify $V_\A^\ast$ with $V_\A$
via the pairing $(x,y)=\ac{[x,y]}$.
Then, the lattice
$V_k\subset V_\A$ becomes self-dual.
For $g=(t,\stwtw abcd)\in G$, set 
$g'=(t^{-1},\frac{1}{ad-bc}\stwtw abcd)$.
Then the above form satisfies
$$[x,y]=[\varrho(g)\cdot x,\varrho(g')\cdot y].$$

Set $dy=dy_0dy_1dy_2$, for $y=(y_0,y_1,y_2)\in V_\A$.
For $\Phi\in{\mathcal S}(V_\adl)$,
we define the Fourier transform $\hat\Phi$ of $\Phi$ by
$$
\hat\Phi(x)=\int_{V_\adl}\Phi(y)(x,y)dy.
$$
Then $\hat\Phi\in{\mathcal S}(V_\adl)$,
and $\hat{\hat\Phi}(x)=\Phi(-x)$.
For $h=\varrho(g)\in H$, set
$\Phi_h(x)=\Phi(h\cdot x)$ and $h'=\varrho(g')$.
Then it is easy to see that the Fourier transform of
$\Phi(h\cdot\quad)$ is $\idn{\chi(h)}^{-3}\hat\Phi(h'\cdot\quad)$.

Set $S_k=V_k\setminus V_k''$, and
$H_\A^0=\{h\in H_\A\mid \idn{\chi(h)}=1\}$.
 For $n\in\Z$ and $\Phi\in\cs(V_\A)$,
we define $\Phi_n(x)=\Phi(\pi^n x)$.
Let
\begin{align}
I^0(\Phi)
&	=
	\int_{H_\adl^0/H_k}
\left(
		\sum_{x\in S_k}\hat\Phi(h'\cdot x)
		-\sum_{x\in S_k}\Phi(h\cdot x)
\right)
	dh,\\
\label{eq:I}
I(\Phi,s)
&	=\frac{1}{2}I^{0}(\Phi)+
			\sum_{n\geq1}q^{-ns}I^0(\Phi_{-n}).
\end{align}

Then by applying the Poisson summation formula to $Z_-(\Phi,s)$,
we obtain the following.
\begin{prop}\label{prop:I}
$$
Z(\Phi,s)=Z_+(\Phi,s)+Z_+(3-s,\hat\Phi)
+I(\Phi,s).
$$
\end{prop}
From now on,
we will study the integral $I^0(\Phi)$
in \S\ref{subsec:ses}--\S\ref{subsec:I^0}
and then compute $I(\Phi,s)$ in \S\ref{subsec:ppf}.

\subsection{The smoothed Eisenstein series}\label{subsec:ses}
To compute $I^0(\Phi)$,
it seems natural to divide the
index set $S_k$ of the summation
into its $H_k$-orbits
and perform integration separately.
However, we cannot put this into practice
because the corresponding integrals diverge.
This is the main difficulty
when one calculates the global zeta functions
of the prehomogeneous vector spaces.
To surmount this problem
Shintani \cite{Shi1} introduced
the smoothed Eisenstein series of $GL(2)$.
Then he determined the principal part
in the case of the space of binary cubic forms
and the space of binary quadratic forms over $\Q$.
Later A.Yukie \cite{yu} generalized
the theory of Eisenstein series
to the groups of products of $GL(n)$'s,
and determined the principal part in some more cases.
In this subsection, we essentially repeat
their argument in our settings.

For $g\in G_\adl$ and $z\in\C$,
we define the Eisenstein series $E(z,g)$ by
\begin{equation*}
E(z,g)=\sum_{\gamma\in G_k/B_k}\idn{\tau(g\gamma)}^{-\frac{z+1}{2}}.
\end{equation*}
This is left $\K Z(G)_\A$-invariant and right $Z(G)_\A G_k$-invariant.
Since the kernel of $\varrho$ is contained in the center of $G$,
we can regard $E(z,g)$ as a function of $h\in H_\adl$.
That is to say, we can write
\begin{equation}
E(z,h)=\sum_{\gamma\in H_k/\B_k}\idn{\tau(h\gamma)}^{-\frac{z+1}{2}}.
\end{equation}
This series converges for ${\rm Re}(z)>1$
and is, in fact, a rational function of $q^{z/2}$.

Write $h=\varrho(\kappa d(t,1)a(\tau)n(u))$.
Then $E(z,h)$ depends only on $\tau,u$, and does not on $\kappa,t$.
Moreover, $E(z,h)$ is right $H_k$-invariant, hence
its value remains unchanged if we replace $u$ by $u+a$ for $a\in k$.
Hence $E(z,h)=E_\tau(z,u)$ has the following Fourier expansion.
\begin{equation}
E_\tau(z,u)=C_0(z,\tau)+\sum_{a\in k^\times}C_a(z,\tau)\ac{au},
\end{equation}
where
$$
C_a(z,\tau)=\int_{\A/k}E_\tau(z,u)\ac{-au}du.
$$

We define
$$
\phi(z)=q^{1-\genus}\frac{\zeta_k(z)}{\zeta_k(z+1)}.
$$
Note that $\phi(z)$ is a rational function of $q^z$ and
is holomorphic in the region $\rez\geq1-\delta$ 
for some $\delta>0$ except for a simple pole
at $z\in\C$ satisfying $q^{1-z}=1$.
The following lemma
about the Fourier coefficients $C_a(z,\tau)$
is well known and we omit the proof.
\begin{lem}
\label{lem:res}
\begin{enumerate}[\rm (i)]
\item
The constant term $C_0(z,\tau)$ has the following explicit formula.
$$C_0(z,\tau)=\idn\tau^{-\frac{z+1}{2}}+
	\idn\tau^{\frac{z-1}{2}}\phi(z).$$
\item
	Let $[\tau]=\sum_{v\in \gM}(\ord_v(\tau))_v$ denote
	the divisor of $\tau$, and $\gc$ a canonical divisor,
	associated with the character $\ac{\cdot}$.
	Then, $C_a(z,\tau)=0$ for all $a\not\in L(\gc-[\tau])$.
	If $a\in L(\gc-[\tau]),a\not=0$, then
$$
C_a(z,\tau)=\idn{\tau}^{\frac{z-1}{2}}\frac{P_a(z,\tau)}{\zeta_k(z+1)},
$$
	where $P_a(z,\tau)$ is a polynomial in $q^{-z}$.
	In particular, $C_a(z,\tau)$ is a holomorphic function
	of $z$ in the half-plane ${\rm Re}(z)>0$.
\end{enumerate}
\end{lem}
Note that the number of $a\in k^\times$ such that
$a\in L(\gc-[\tau])$ is finite.
From this, we immediately obtain the following.

\begin{cor}\label{cor:rho0}
The function $E(z,\tau)$ is holomorphic in the region $\rez>0$
with an exception of a simple pole
at $z\in\C$ satisfying $q^{1-z}=1$, and
\begin{equation*}
\lim_{z\to1}(1-q^{1-z})E(z,\tau)
=	\lim_{z\to1}(1-q^{1-z})\phi(z)
=	\frac{q^{2-2\genus}\gC_k}{\zeta_k(2)}.
\end{equation*}
\end{cor}
We denote the value in the formula of Corollary \ref{cor:rho0} by $\rho_0$.
Let $\psi$ be an entire function such that for any $c_1,c_2\in \R$ and $N>0$,
$$
\sup_{c_1<{\rm Re}(w)<c_2}\left(1+|w|^N\right)|\psi(w)|<\infty.
$$
Let ${\rm Re}(w)>1$.
Following \cite{Shi1}, we define the smoothed Eisenstein series
$\E(w,h)$ by	
\begin{equation}
\E(w,h)=\frac{1}{2\pi i}\int_{{\rm Re}(z)=r}
\frac{E(z,h)}{w-z}\psi(z)dz
\end{equation}
for some $1<r<{\rm Re}(w)$.
Note that this integral does not depend on the choice of $r$.
Similarly as $E(z,h)$,
$\E(w,h)$ has the Fourier expansion
$$
\E(w,h)=\sum_{a\in k}\Cs_a(w,\tau)\ac{au},
$$
where
$$
\Cs_a(w,\tau)=
	\frac{1}{2\pi i}\int_{{\rm Re}(z)=r}
	\frac{C_a(z,\tau)}{w-z}\psi(z)dz.
$$

This $\E(w,h)$ satisfies the following property.

\begin{lem}\label{lem:ses}
\begin{enumerate}[\rm (i)]
\item
	As a function of $w$,
 	$\E(w,h)$ is holomorphic
 	in the region ${\rm Re}(w)>1$.
\item
	For any $w$ such that $\rew>1$,
	$\E(w,h)\ll\idn{\tau(h)}^{(\rew-1)/2}$.
\item\label{it:notice}
 	$\lim_{w\to1+0}(1-q^{1-w})\E(w,h)=\rho_0\psi(1)$.
\item
	For any $M>1$,
	$$
	\underset{h\in H_\A}{\sup_{1<\rew<M}}
		|(1-q^{1-w})\E(w,h)|<\infty.
	$$
\end{enumerate}
\end{lem}
Since the proof is similar to that
of Lemma 3.2 of \cite{Dat2},
we omit it.
(Note that the convergence of (\ref{it:notice}) is not uniform.)
As a result of this, we have the following.
\begin{cor}\label{cor:exchg}
	For $f\in L^1(H_\adl^0/H_k)$,
	$$
		\lim_{w\to1+0}(1-q^{1-w})
		\int_{H_\adl^0/H_k}f(h)\E(w,h)dh
		=\rho_0\psi(1)
		\int_{H_\adl^0/H_k}f(h)dh.
	$$
\end{cor}

\subsection{Decomposition of $I^0(\Phi,w)$}
Set
\begin{equation}
I^0(\Phi,w)=
	\int_{H_\adl^0/H_k}\sum_{x\in S_k}
	\left(\hat\Phi(h'\cdot x)-\Phi(h\cdot x)\right)
	\E(w,h)dh.
\end{equation}
By Corollary \ref{cor:exchg}, we have
\begin{equation*}
	\lim_{w\to1+0}(1-q^{1-w})I^0(\Phi,w)=
	\rho_0\psi(1)I^0(\Phi).
\end{equation*}
We have the following lemma on the structure of $S_k$.
This can be easily proved and
we simply state the result here.

\begin{lem}\label{lem:S-st}
Let $S_k^i,i=0,1,2$, be the subsets of $V_k$ given by
$$
S_k^0=\{0\},\quad
S_k^1=\{x\in V_k\mid x\not=0,P(x)=0\},\quad
S_k^2=\{x\in V_k\mid k(x)=k,P(x)\not=0\}.$$
Then $S_k=S_k^0\amalg S_k^1\amalg S_k^2$ and moreover,
$$
S^1_k=H_k\times_{B_k}\{(0,0,a)|a\in k^\times\},\qquad
S^2_k=H_k\times_{\widehat{\T}_k}\{(0,a,0)|a\in k^\times\}.
$$
\end{lem}

\begin{defn}\label{def:J}
For $i=0,1,2$, we define 
\begin{equation*}
J_i(\Phi,w)=
	\int_{H_\adl^0/H_k}\sum_{x\in S_k^i}\Phi(h\cdot x)\E(w,h)dh.
\end{equation*}
\end{defn}
The next lemma shows that each $J_i(\Phi,w)$
converges and is holomorphic in the right half-plane $\rew>1$.
\begin{lem}
The integral
$$	\int_{H_\adl^0/H_k}
		\sum_{x\in V_k}|\Phi(h\cdot x)|
		|\E(w,h)|dh
$$
converges absolutely and locally uniformly in the region ${\rm Re}(w)>1$.
\end{lem}
\begin{proof}
Set $\Theta^0=\{g\in\Theta|\idn{\det(\varrho(g))}=1\}$.
For $\idn{t\tau}=1$ and $\idn{\tau}\leq q^{2\genus}$,
by Lemma \ref{lem:weil},
\begin{align*}
	\sum_{x\in V_k}|\Phi(\varrho(d(t,1)a(\tau))x)|
&		=\sum_{x\in V_k}|\Phi(tx_0,t\tau x_1,t\tau^2 x_2)|\\
&	\ll\max\{1,\idn{t}^{-1}\}\max\{1,\idn{t\tau}^{-1}\}
		\max\{1,\idn{t\tau^2}^{-1}\}
		\ll\idn{\tau}^{-1}.
\end{align*}
Then, by Lemma \ref{lem:ses}(ii),
the integral is bounded by a constant multiple of
\begin{align*}
 \int_{\K}d\kappa
		\underset{\idn{t\tau}=1,\idn{\tau}\leq q^{2\genus}}
			{\idlint{t}\idlint{\tau}}
		\adlint{u}\ \ \idn{\tau}^{\frac{{\rm Re}(w)-1}{2}}
&=	\underset{\idn{\tau}\leq q^{2\genus}}
			{\int_{\idl/k^\times}}
		\idn{\tau}^{\frac{{\rm Re}(w)-1}{2}}d^\times\tau
<	\infty.
\end{align*}
\end{proof}

Note that if $h\in H_\A^0$, $h=h'$.
Then, by Lemma \ref{lem:S-st},
we have the following.
\begin{prop}\label{prop:Jsum}
\begin{equation}\label{eq:Jsum}
I^0(\Phi,w)
=	J_0(\hat\Phi,w)-J_0(\Phi,w)
+	J_1(\hat\Phi,w)-J_1(\Phi,w)
+	J_2(\hat\Phi,w)-J_2(\Phi,w).
\end{equation}
\end{prop}
From now on,
we will compute $J_0,J_1,J_2$ in
\S\ref{subsec:j0},\ref{subsec:j1},\ref{subsec:j2},
respectively.
For this purpose, we will introduce some notation.
For meromorphic functions $f_1(w),f_2(w)$,
we will use the notation $f_1\sim f_2$
if $f_1-f_2$ can be continued meromorphically to 
a right half plane $\rew>1-\sigma$
for some $\sigma>0$ 
and becomes holomorphic in the region.

\subsection{Computation of $J_0(\Phi,w)$}
\label{subsec:j0}

\begin{prop}\label{prop:j0}
\begin{equation*}
J_0(\Phi,w)\sim
\Phi(0)\frac{\psi(w)}{1-q^{\frac{1-w}{2}}}.
\end{equation*}
\end{prop}

\begin{proof}
By the formula of \S\ref{subsec:dh1},\eqref{eq:edh0},
we have
\begin{align*}
J_0(\Phi,w)
&=	\Phi(0)\int_{H_\adl^0/H_k}
	\left(\cmpint
	\frac{\sum_{\gamma\in H_k/\B_k}\idn{\tau(h\gamma)}^{-\frac{z+1}{2}}}{w-z}
	\psi(z)dz\right)dh\\
&=	\Phi(0)\int_{H_\adl^0/\B_k}
	\left(\cmpint
	\frac{\idn{\tau(h)}^{-\frac{z+1}{2}}}{w-z}
	\psi(z)dz\right)dh\\
&=	\Phi(0)
	\underset{\idn{t\tau}=1}{\idlint{t}\idlint{\tau}}\adlint{u}
	\left(\cmpint
	\frac{\idn{\tau}^{-\frac{z+1}{2}}}{w-z}
	\psi(z)dz\right)\idn{\tau}.
\intertext
{By replacing $\tau$ by $\tau/t$,}
&= 	\Phi(0)
	\idlint{t}
	\left(\cmpint
	\frac{\idn{t}^{\frac{z-1}{2}}}{w-z}
	\psi(z)dz\right)\\
&=	\Phi(0)
	\left(\int_{\adl^+/k^\times}\hspace*{\idllength}d^\times t
	+\int_{\adl^-/k^\times}\hspace*{\idllength}d^\times t\right)
	\left(\cmpint
	\frac{\idn{t}^{\frac{z-1}{2}}}{w-z}
	\psi(z)dz\right),
\end{align*}
where,
$\adl^+=\{t\in\idl|\idn{t}>1\},\adl^-=\{t\in\idl|\idn{t}\leq1\}$.
Let $r_1<1$. The former integral is equal to
\begin{align*}
	\int_{\adl^+/k^\times}\hspace*{\idllength}d^\times t
	\ \ \frac{1}{2\pi i}
		\int_{\rez=r_1}\frac{\idn{t}^{\frac{z-1}{2}}}{w-z}\psi(z)dz
&=	\sum_{m=1}^\infty\frac{1}{2\pi i}
		\int_{\rez=r_1}\frac{q^{\frac{z-1}{2}m}}{w-z}\psi(z)dz\\
&=	\frac{1}{2\pi i}\int_{\rez=r_1}
		\frac{(q^{\frac{1-z}{2}}-1)^{-1}}{w-z}\psi(z)dz
\end{align*}
and is holomorphic for $\rew>r_1$. Hence,
\begin{align*}
J_0(\Phi,w)
&\sim
	\Phi(0)
	\int_{\adl^-/k^\times}\hspace*{\idllength}d^\times t
	\ \ \frac{1}{2\pi i}\int_{\rez=r}
	\frac{\idn{t}^{\frac{z-1}{2}}}{w-z}\psi(z)dz\\
&=	\Phi(0)
	\ \ \frac{1}{2\pi i}\int_{\rez=r}
	\frac{(1-q^{\frac{1-z}{2}})^{-1}}{w-z}\psi(z)dz\\
&=	\Phi(0)\frac{\psi(w)}{1-q^{\frac{1-w}{2}}}
+	\Phi(0)
	\ \ \frac{1}{2\pi i}\int_{{\rm Re}(z)=r_2>{\rm Re}(w)}
	\frac{(1-q^{\frac{1-z}{2}})^{-1}}{w-z}\psi(z)dz\\
&\sim
	\Phi(0)\frac{\psi(w)}{1-q^{\frac{1-w}{2}}}.
\end{align*}
We have thus proved the proposition.
\end{proof}

%

\subsection{Computation of $J_1(\Phi,w)$}
\label{subsec:j1}
For $\Phi\in\mathcal{S}(V_\adl)$, we define 
$M\Phi\in\mathcal{S}(V_\A)$ by
\begin{equation}
M\Phi(x)=
	\int_\K\Phi(\varrho(\kappa)\cdot x)d\kappa.
\end{equation}
Since $Z(\Phi,s)=Z(M\Phi,s)$,
we assume that $\Phi=M\Phi$ for the rest of this section.

For $\Psi\in\mathcal{S}(\adl)$ and $s\in\C$,
Tate's zeta function $\Sigma(\Psi,s)$ is defined by
\begin{equation}
\Sigma(\Psi,s)=\int_{\idl}\idn{t}^s\Psi(t)d^\times t.
\end{equation}
It is well known (see \cite{Ta},\cite{We}) that
$\Sigma(\Psi,s)$ can be written as follows:
\begin{equation}
\Sigma(\Psi,s)=P(\Psi,s)+
\frac{\hat\Psi(0)}{1-q^{1-s}}-\frac{\Psi(0)}{1-q^{-s}},
\end{equation}
where $P(\Psi,s)$ is a polynomial in $q^s,q^{-s}$.

For $\Phi\in\mathcal{S}(V_\adl)$,
define $R_1\Phi\in\mathcal{S}(\adl)$ as
$R_1\Phi(t)=\Phi(0,0,t)$.

\begin{prop}\label{prop:j1}
\begin{equation*}
	J_1(\Phi,w)
	\sim
		\phi(w)\psi(w)
		\Sigma\left(R_1\Phi,\frac{w+1}{2}\right).
\end{equation*}
\end{prop}

\begin{proof}
By Lemma \ref{lem:S-st} and the formula
\S\ref{subsec:dh1},\eqref{eq:edh1}, we have
\begin{align*}
J_1(\Phi,w) 
&=	\int_{H_\adl^{0}/H_k}
	\sum_{\gamma\in H_k/\B_k}\sum_{a\in k^\times}
	\Phi(h\gamma\cdot(0,0,a))\E(w,h)dh\\
&=	\int_{H_\adl^{0}/\B_k}
	\sum_{a\in k^\times}\Phi(h\cdot(0,0,a))\E(w,h)dh\\
&=	\int_\K d\kappa
	\underset{\idn{t\tau}=1}{\idlint{t}\idlint{\tau}}\adlint{u}
	\sum_{a\in k^\times}\Phi(\varrho(\kappa)\cdot(0,0,at\tau^2))
	\E(w,a(\tau)n(u))\idn{\tau}\\
&=	\underset{\idn{t\tau}=1}{\int_{\idl}d^\times t\idlint{\tau}}
	\ \ R_1\Phi(t\tau^2)
	\Cs_0(w,\tau)\idn{\tau}.
\end{align*}
Here, in the last expression,
we replaced $\E(w,h)$ by $\Cs_0(w,\tau)$
by the orthogonality of characters.
By replacing $\tau$ by $\tau/t$ and afterwards $t$ by $\tau^2/t$,
we have
\begin{align*}
J_1(\Phi,w) 
&=	\int_{\idl}d^\times t\underset{\idn{\tau}=1}{\idlint{\tau}}
	\ \ R_1\Phi(t)
	\Cs_0\left(w,\frac{t}{\tau}\right)\idn{\frac{t}{\tau}}\\
&=	\int_{\idl}d^\times t\ \ R_1\Phi(t)
	\frac{1}{2\pi i}\int_{{\rm Re}(z)=r}
	\frac{\idn{t}^{-\frac{z-1}{2}}+\idn{t}^{\frac{z+1}{2}}\phi(z)}{w-z}
	\psi(z)dz.
\end{align*}
We now break the above integral into two parts.
Let $r_1<-1$.
The first part is equal to
$$	
	\frac{1}{2\pi i}\int_{{\rm Re}(z)=r_1}
		\Sigma\left(R_1\Phi,-\frac{z-1}{2}\right)
	\frac{\psi(z)}{w-z}dz.
$$
By the theory of Tate's zeta function,
$\Sigma\left(R_1\Phi,-\frac{z-1}{2}\right)$
is a holomorphic function in the region $\rez<-1$.
Therefore this part is an entire function of $w$. Hence,
\begin{align*}
J_1(\Phi,w) 
&\sim	
	\frac{1}{2\pi i}\int_{\rez=r}
	\Sigma\left(R_1\Phi,\frac{z+1}{2}\right)
	\frac{\phi(z)}{w-z}
	\psi(z)dz\\
&\sim	
	\phi(w)\psi(w)
	\Sigma\left(R_1\Phi,\frac{w+1}{2}\right).
\end{align*}
\end{proof}

\subsection{Unstable distributions}
\label{subsec:ud}
Here, we introduce some distributions
and consider its analytic properties,
which will be needed later.
The results in this subsection are
essentially due to \cite{yu1} Section 2.
\begin{defn}
For $s,w\in\C$ and $\Psi\in{\mathcal S}({\rm Aff}_\A^2)$, 
we define
\begin{align*}
T(\Psi,s,w)&=	\int_{\idl}\int_{\adl}
					\idn{t}^s\alpha(u)^{-w}\Psi(t,tu)
				dud^\times t,\\
T^+(\Psi,s,w)&=	\int_{\idl}\int_{\adl}
					\lambda_+(t)\idn{t}^s\alpha(u)^{-w}\Psi(t,tu)
				dud^\times t,\\
T^-(\Psi,s,w)&=	\int_{\idl}\int_{\adl}
					\lambda_-(t)\idn{t}^s\alpha(u)^{-w} \Psi(t,tu)
				dud^\times t,\\
T^0(\Psi,w)&=	\int_{\adl^0}\int_{\adl}
					\alpha(u)^{-w} \Psi(t,tu)
				dud^\times t,
\end{align*}
where $\alpha(u)$ is
introduced in \S\ref{subsec:dh1}.
\end{defn}
\begin{lem}\label{lem:T+rat}
\begin{enumerate}[\rm(i)]
\item\label{it:t1}
$T^+(\Psi,s,w)$ converges absolutely and locally uniformly
for all $s,w\in\C$,
and
$T^-(\Psi,s,w)$ converges absolutely and locally uniformly
for all $\res+\rew>2, \res>2$.
In particular,
$T^0(\Psi,w)$ converges absolutely and locally uniformly
for all $w\in\C$.
\item\label{it:t2}
As a function of $s$,
$T^+(\Phi,s,w)$ is a polynomial in $q^s$.
\end{enumerate}
\end{lem}
\begin{proof}
Let $\sigma=\res, \sigma_1=\rew$.
Let $f,g\geq0$ be
Schwartz-Bruhat functions on $\A$
such that
$|\Psi(x_1,x_2)|\leq f(x_1)g(x_2)$
for $x_1,x_2\in\A$.
By changing $u$ to $ut^{-1}$,
we have
\begin{equation*}
T^-(\Psi,s,w)
\ll
\idlint{t}\int_\A du\ \ 
	\lambda_-(t)\idn{t}^{\sigma-1}
	\sum_{x\in k^\times}
	f(tx)g(u)
	\alpha(t^{-1}x^{-1}u)^{-\sigma_1}.
\end{equation*}
Also, the same argument as the proof of Lemma \ref{lem:t}
shows that there exists an integer $N>0$ such that
\begin{equation*}
T^+(\Psi,s,w)
\ll
\underset{\idn t\leq q^N}{\idlint{t}}
\int_\A du\ \ 
	\lambda_+(t)\idn{t}^{\sigma-1}
	\sum_{x\in k^\times}
	f(tx)g(u)
	\alpha(t^{-1}x^{-1}u)^{-\sigma_1}.
\end{equation*}
Let $C_1=\supp f, C_2=\supp g$.
We will give an estimate
$\alpha(t^{-1}x^{-1}u)^{-\sigma_1}$
for $t\in\A^\times, x\in k^\times, u\in\A$
such that $tx\in C_1, u\in C_2$.
Since $\alpha(\cdot)\leq 1$,
we have
$\alpha(t^{-1}x^{-1}u)^{-\sigma_1}\leq 1$
if $\sigma_1\geq 0$.
Let $\sigma_1\leq 0$.
By the definition of $\alpha(\cdot)$, we have
\begin{align*}
\alpha(t^{-1}x^{-1}u)
&=	\prod_{v\in\gM}\sup(1,|t^{-1}x^{-1}u|_v)\\
&=	\prod_{v\in\gM}|t^{-1}x^{-1}|_v\sup(|tx|_v,|u|_v)
=	\idn{t}^{-1}\prod_{v\in\gM}\sup(|tx|_v,|u|_v).
\end{align*}
Since $C_1,C_2$ are compact, 
$\prod_{v\in\gM}\sup(|tx|_v,|u|_v)$
is bounded by a constant. Hence,
\begin{equation*}
\alpha(t^{-1}x^{-1}u)^{-\sigma_1}\ll \idn{t}^{\sigma_1}
\end{equation*}
Therefore (\ref{it:t1}) follows from Lemma \ref{lem:weil}.

Let $\Psi_n(\cdot)=\Psi(\pi^n\cdot)$.
Then we have
\begin{equation*}
T^+(\Psi,s,w)
	=\frac{1}{2}T^0(\Psi,w)
	+\sum_{n=1}^{N}q^{ns}T^0(\Psi_n,w).
\end{equation*}
This establishes (\ref{it:t2}).
\end{proof}

From now on,
we will give the explicit formula of $T(\Psi,s,w)$.
Since $\alpha(u)$ is a product of local factors $\alpha_v(u)$,
we can also define local distributions $T_v$ for 
$\Psi_v\in{\mathcal S}(k_v^2)$:
\begin{equation*}
T_v(\Psi_v,s,w)
	=\int_{k_v^\times}\int_{k_v}|t_v|_v^s
		\alpha_v(u_v)^{-w}
		\Psi_v(t_v,t_vu_v)
					\ du_v d^\times t_v.
\end{equation*}
Then by \eqref{productmeasures},
$$
T(\Psi,s,w)=
q^{1-\genus}\gC_k^{-1}\prod_{v\in\gM}T_v(\Psi_v,s,w)
$$
for $\Psi=\prod_{v\in\gM}\Psi_v$.
We recall Proposition 2.8 and 2.9 of \cite{yu1}.

\begin{lem}
\begin{enumerate}[\rm (i)]
\item
	$T_v(\Psi_v,s,w)$ converges absolutely and locally uniformly
	in the region $\res+\rew>1,
	\res>0$, and is holomorphic in the region.
\item
	If $\Psi_v$ is the characteristic function of $\co_v^2$,
	\begin{equation*}
		T_v(\Psi_v,s,w)
		=\frac{1-q_v^{-(s+w)}}{(1-q_v^{-s})(1-q_v^{-(s+w-1)})}.
	\end{equation*}
\item
	For any $\Psi_v$,
	$(1-q_v^{-s})(1-q_v^{-(s+w-1)})T_v(\Psi_v,s,w)$
	becomes a polynomial in $q_v^{\pm s},q_v^{\pm w}$.
\end{enumerate}
\end{lem}

This lemma implies

\begin{lem}\label{lem:Trat}
$T(\Psi,s,w)$ is holomorphic in the region 
$\res+\rew>2,\res>1$,
and moreover, is a rational function of $q^{-s},q^{-w}$.
More precisely, suppose that $\Psi=\otimes\Psi_v$
and that $S$ a finite set of places such that 
$\Psi_v$ are the characteristic functions of $\co_v^2$ for $v\not\in S$,
then
\begin{equation*}
T(\Psi,s,w)=T_S(\Psi,s,w)
	\frac{\zeta_{k,S}(s)\zeta_{k,S}(s+w-1)}{\zeta_{k,S}(s+w)},
\end{equation*}
where
$T_S(\Psi,s,w)=
q^{1-\genus}\gC_k^{-1}\prod_{v\in S}T_v(\Psi_v,s,w)$ and
\begin{equation*}
\zeta_{k,S}(s)=
	\prod_{v\not\in S}(1-q_v^{-s})^{-1}
\end{equation*}
is the truncated Dedekind zeta function.
\end{lem}


Let us define the distribution
$\tilde{T}(\Psi,s)$ by
\begin{equation*}
\tilde{T}(\Psi,s)=\frac{1}{\log q}\frac{d}{dw}T(\Psi,s,w)|_{w=0}.
\end{equation*}
Also we define $\tilde{T}^+(\Psi,s),\tilde{T}^-(\Psi,s)$ in similar way.
These are rational functions of $q^s$.
For later purposes, we will state the pole structure of
$\tilde{T}(\Psi,s)$. Since
\begin{multline*}
\tilde{T}(\Psi,s)
=	\zeta_{k,S}(s-1)\frac{d}{dw}T_S(\Psi,s,w)|_{w=0}
+	T_S(\Psi,s,0)
		\left(\zeta'_{k,S}(s-1)-
			\zeta_{k,S}(s-1)\frac{\zeta'_{k,S}(s)}{\zeta_{k,S}(s)}
		\right)
\end{multline*}
and hence,
\begin{lem}\label{lem:pd}
$\tilde{T}(\Psi,s)$
is a rational function of $q^s$,
and holomorphic in the region $\res>2$.
It has at most double pole at $s=2$,
and $(1-q^{2-s})^2\tilde{T}(\Psi,s)$
is holomorphic in the region $\res>1$.
\end{lem}

\subsection{Computation of $J_2(\Phi,w)$}
\label{subsec:j2}

For $\Phi\in{\mathcal S}(V_\adl)$,
we define $R_2\Phi\in{\mathcal S}({\rm Aff}_\A^2)$ by
$R_2\Phi(t,u)=\Phi(0,t,u)$.

\begin{prop}\label{prop:j2}
\begin{equation*}
J_2(\Phi,w)
\sim	\phi(w)\psi(w)
		\left(\frac{1}{1-q^{\frac{1-w}{2}}}-\frac{1}{2}\right)
		T^0\left(R_2\Phi,\frac{1-w}{2}\right).
\end{equation*}
\end{prop}

\proof
By the formula of \S\ref{subsec:dh1},\eqref{eq:edh2},
\begin{align*}
J_2(\Phi,w) 
&=	\int_{H_\adl^{0}/H_k}
	\sum_{\gamma\in H_k/\widehat{\T}_k}\sum_{a\in k^\times}
	\Phi(h\gamma\cdot(0,a,0))\E(w,h)dh\\
&=	\int_{H_\adl^{0}/\widehat{\T}_k}
	\sum_{a\in k^\times}\Phi(h\cdot(0,a,0))\E(w,h)dh\\
&=	\int_\K d\kappa
	\int_{\adl}du
	\int d^\times\tau
	\underset{\idn{t\tau}=1}{\idlint{t}}
	\sum_{a\in k^\times}\Phi(\varrho(\kappa)\cdot(0,at\tau,aut\tau))
	\E(w,n(u)a(\tau))\\
&=	\int_{\adl}du
	\int d^\times\tau
	\int_{\adl^0}d^\times t\ \ 
	R_2\Phi(t,tu)\E(w,a(\tau)n(\frac{u}{\tau})),
\end{align*}
where the last transformation can be obtained
by changing $t$ to $t/\tau$, and including
the sum $\sum_{a\in k^\times}$ into the integration of $t$.
Here, we used the notation 
$$
\int d^\times\tau
=		\underset{\idn{\tau}<\alpha(u)}{\idlint{\tau}}
		+\frac{1}{2}\underset{\idn{\tau}=\alpha(u)}{\idlint{\tau}}
$$
for simplicity. 
Now, we show the following claim.

\begin{clm}
Set
$$
\E'(w,h)=\sum_{a\in k^\times}\Cs_a(w,\tau)\ac{au}
=\E(w,h)-\Cs_0(w,\tau).
$$
Then
$$
P_1=
\int_{\adl}du
\int d^\times\tau
\int_{\adl^0}d^\times t\ \ 
R_2\Phi(t,tu)\E'(w,a(\tau)n(\frac{u}{\tau}))
$$
is a holomorphic function of $w$ in the half-plane $\rew>0$.
\end{clm}

\noindent
{\it Proof of Claim.}
By the Fourier expansion of the Eisenstein series,
$$
P_1=\sum_{a\in k^\times}
\int_{\adl}du
\int d^\times\tau
\int_{\adl^0}d^\times t\ \ 
R_2\Phi(t,tu)\Cs_a(w,\tau)
\ac{\frac{ua}{\tau}}.
$$
We can see that there are finitely many $a$'s such that
$\Cs_a(w,\tau)\not=0$ for some $\tau$ with
$\idn{\tau}<\alpha(u)$
for the following reason.
\begin{enumerate}[\rm (i)]
\item
	For each $\tau$, the number of $a$'s such that
	$\Cs_a(w,\tau)\not=0$ is finite.
	Moreover, for the same $k^\times\idl(\emptyset)$-coset in $\A^\times$
	the Fourier coefficients $\Cs_a(w,\tau)$ are equal.
\item\label{item:cpt}
	By Lemma \ref{lem:res}, $\Cs_a(w,\tau)=0$ for all $a\in k^\times$
	when $\idn{\tau}<q^{2-2\genus}$.
\item
	The set $\{\tau\in\idl/k^\times\idl(\emptyset)\mid
	q^{2-2\genus}\leq\idn{\tau}\leq\alpha(u)\}$
	is a finite set.
\end{enumerate}
Therefore, it is enough to prove that
each integral in the sum $\sum_{a\in k^\times}$
is holomorphic in the region $\rew>0$.
On the other hand, (\ref{item:cpt}) implies that
the region of integration with respect to $\tau$ is compact.
Since $\Cs_a(w,\tau)$ is holomorphic in the region $\rew>0$,
 we obtain the claim.
\hspace*{\fill}{$\square$}

Now by the claim above, we have
\begin{align*}
J_2(\Phi,w)
&=	\int_{\adl}du\int_{\adl^0}d^\times t\int d^\times\tau\ \ 
	R_2\Phi(t,tu)
	\frac{1}{2\pi i}\int_{\rez=r}
	\frac{\idn{\tau}^{-\frac{z+1}{2}}+
	\idn{\tau}^{\frac{z-1}{2}}\phi(z)}{w-z}
	\psi(z)dz.
\end{align*}

We break the above integral into two parts.
Note that
\begin{equation*}
\int d^\times\tau\ \idn{\tau}^z
=	\left(\frac{1}{1-q^{-z}}-\frac{1}{2}\right)\alpha(u)^z
\end{equation*}
for $\rez>0$.
Let $r_1<-1$.
Then the first one is equal to
\begin{multline*}
\int_{\adl}du\int_{\adl^0}d^\times t\ \ 
	R_2\Phi(t,tu)
	\int d^\times\tau
	\ \ \frac{1}{2\pi i}\int_{\rez=r_1}
	\idn{\tau}^{-\frac{z+1}{2}}
	\frac{\psi(z)}{w-z}dz\\
=	\frac{1}{2\pi i}\int_{\rez=r_1}
	\left(\frac{1}{1-q^{\frac{z+1}{2}}}-\frac{1}{2}\right)
	T^0\left(R_2\Phi,\frac{z+1}{2}\right)
	\frac{\psi(z)}{w-z}
	dz
\end{multline*}
and hence, is holomorphic in the region $\rew>r_1$.
Therefore,
\begin{align*}
J_2(\Phi,w)
&\sim
	\int_{\adl}du\int_{\adl^0}d^\times t\int d^\times\tau\ \ 
	R_2\Phi(t,tu)
	\frac{1}{2\pi i}\int_{\rez=r}
	\frac{\idn{\tau}^{\frac{z-1}{2}}\phi(z)}{w-z}
	\psi(z)dz\\
&=	\phi(w)	\psi(w)
	\int_{\adl}du\int_{\adl^0}d^\times t\ \ 
	R_2\Phi(t,tu)
	\int d^\times\tau
	\idn{\tau}^{\frac{w-1}{2}}\\
&\quad\
+	\int_{\adl}du\int_{\adl^0}d^\times t\int d^\times\tau\ 
	R_2\Phi(t,tu)
	\frac{1}{2\pi i}\int_{{\rm Re}(z)=r_2>{\rm Re}(w)}
	\frac{\idn{\tau}^{\frac{z-1}{2}}\phi(z)}{w-z}
	\psi(z)dz\\
&=	\phi(w)\psi(w)
	\left(\frac{1}{1-q^{\frac{1-w}{2}}}-\frac{1}{2}\right)
	T^0\left(R_2\Phi,\frac{1-w}{2}\right)\\
&\quad
+	\int_{\adl}du\int_{\adl^0}d^\times t\int d^\times\tau\ 
	R_2\Phi(t,tu)
	\frac{1}{2\pi i}\int_{{\rm Re}(z)=r_2>{\rm Re}(w)}
	\frac{\idn{\tau}^{\frac{z-1}{2}}\phi(z)}{w-z}
	\psi(z)dz.
\end{align*}
Now, Similarly to $P_2$,
we can show that 
the second part of the last
expression is an entire function of $w$.
Hence we have the proposition.
\hspace*{\fill}$\square$

\subsection{Explicit evaluation of $I^0(\Phi)$}
\label{subsec:I^0}

Now, we turn to $I^0(\Phi)$.
Recall that 
\begin{equation}\label{eq:I^0}
	\lim_{w\to1+0}(1-q^{1-w})I^0(\Phi,w)=
	\rho_0\psi(1)I^0(\Phi)
\end{equation}
by Corollary \ref{cor:exchg}.
We will compute $I^0(\Phi)$ by using
Proposition \ref{prop:Jsum} and \ref{prop:j0}, \ref{prop:j1}, \ref{prop:j2}.
Recall that $\Sigma(\Psi_1,w)$ has a simple pole at $w=1$,
and $T^0(\Psi_2,w)$ is holomorphic at $w=0$.
We will write their Laurent expansion
in $q^{1-w}, q^{-w}$, respectively,
at the values $q^{1-w}=1,q^{-w}=1$, by
\begin{align*}
\Sigma(\Psi_1,w)=	\sum_{i=-1}^\infty\Sigma_{(i)}(\Psi_1,1)(1-q^{1-w})^i,
\qquad
T^0(\Psi_2,w)	=	\sum_{i=0}^\infty T_{(i)}^0(\Psi_2,0)(1-q^{-w})^i.
\end{align*}
Obviously,
	$T_{(0)}^0(\Psi_2,0)=T^0(\Phi_2,0),
	T_{(1)}^0(\Psi_2,0)=\frac{1}{\log q}\frac{d}{dw}
		T^0(\Psi_2,w)|_{w=0}$.

Since the limit \eqref{eq:I^0} exists,
the double poles of the right-hand side of \eqref{eq:Jsum}
at $w=1$ cancel out on the whole.
Hence we get the following lemma.
Note that some straightforward calculation also shows the following equality.
(See \cite{Dat1}.)
\begin{lem}
\begin{equation*}
	T_{(0)}^0(R_2\hat\Phi,0)-T_{(0)}^0(R_2\Phi,0)
		=\Sigma_{(-1)}(R_1\Phi,1)-\Sigma_{(-1)}(R_1\hat\Phi,1).
\end{equation*}
\end{lem}

Hence, we have
\begin{align*}
\begin{split}
I(\Phi,w)
&\sim	\frac{\psi(w)}{1-q^{\frac{1-w}{2}}}
		(\hat\Phi(0)-\Phi(0))\\
&\quad
+	\phi(w)\psi(w)
	\left\{\frac{\Sigma_{(-1)}(R_1\hat\Phi,1)}{2}+\Sigma_{(0)}(R_1\hat\Phi,1)
-\frac{\Sigma_{(-1)}(R_1\Phi,1)}{2}-\Sigma_{(0)}(R_1\Phi,1)\right\}\\
&\quad
-	\phi(w)\psi(w)
	\left\{
		T_{(1)}^0(R_2\hat\Phi,0)-T_{(1)}^0(R_2\Phi,0)
	\right\}.
\end{split}
\end{align*}
Then, together with \eqref{eq:I^0}, we obtain the following.

\begin{prop}\label{prop:I^0}
\begin{align*}
I^0(\Phi)
&=	\frac{2}{\rho_0}(\hat\Phi(0)-\Phi(0))\\
&\quad
+	\left\{
		\frac{\Sigma_{(-1)}(R_1\hat\Phi,1)}{2}+\Sigma_{(0)}(R_1\hat\Phi,1)
		-\frac{\Sigma_{(-1)}(R_1\Phi,1)}{2}-\Sigma_{(0)}(R_1\Phi,1)
	\right\}\\
&\quad
-	T_{(1)}^0(R_2\hat\Phi,0)+T_{(1)}^0(R_2\Phi,0).
\end{align*}
\end{prop}

\subsection{The principal part formula}\label{subsec:ppf}

\begin{thm}\label{th:gz}
Suppose that $\Phi=M\Phi$.
Then $Z(\Phi,s)$ is a rational function of $q^s$.
More precisely,
\begin{equation*}
	Z(\Phi,s)=Z_+(\Phi,s)+Z_+(\hat\Phi,3-s)
	+I(\Phi,s),
\end{equation*}
where 
$Z_+(\Phi,s)$ and $Z_+(\hat\Phi,3-s)$ are
polynomials in $q^{s}$ and $q^{-s}$, respectively,
and $I(\Phi,s)$ is given by
\begin{equation}\label{eq:II}
\begin{split}
&	\frac{2}{\rho_0}
	\left\{
		\left(\frac{1}{2}+\frac{q^{3-s}}{1-q^{3-s}}\right)\hat\Phi(0)
		-\left(\frac{1}{2}+\frac{q^{-s}}{1-q^{-s}}\right)\Phi(0)
	\right\}\\
&\quad
+	
	\left\{
		\left(\frac{\Sigma_{(-1)}(R_1\hat\Phi,1)}{2}+\Sigma_{(0)}(R_1\hat\Phi,1)\right)
		\left(\frac{1}{2}+\frac{q^{2-s}}{1-q^{2-s}}\right)
		-\Sigma_{(-1)}(R_1\hat\Phi,1)\frac{q^{2-s}}{(1-q^{2-s})^2}
	\right\}\\
&\quad
-	
	\left\{
		\left(\frac{\Sigma_{(-1)}(R_1\Phi,1)}{2}+\Sigma_{(0)}(R_1\Phi,1)\right)
		\left(\frac{1}{2}+\frac{q^{1-s}}{1-q^{1-s}}\right)
		+\Sigma_{(-1)}(R_1\Phi,1)\frac{q^{1-s}}{(1-q^{1-s})^2}
	\right\}\\
&\quad
-	(
		\tilde T^+(R_2\hat\Phi,3-s)
		-\tilde T^-(R_2\Phi,s)
	).
\end{split}
\end{equation}
\end{thm}

\begin{proof}
We will compute $I(\Phi,s)$
by
\eqref{eq:I} and 
Proposition \ref{prop:I^0}.
Note that $\widehat{\Phi_n}(t)=q^{-3n}\hat\Phi_{-n}(t)$.
Then we have
\begin{equation}\label{eq:calc1}
	\Phi_n(0)=\Phi(0),\qquad \widehat{\Phi_n}(0)=q^{-3n}\hat\Phi(0).
\end{equation}
Also since
\begin{equation*}
	\Sigma(\Phi_n,w)=q^{-nw}\Sigma(\Phi,w),\qquad
	\Sigma(R_1\widehat{\Phi_n},w)=q^{-n(3-w)}\Sigma(R_1\hat\Phi,w),
\end{equation*}
we have
\begin{align}\label{eq:calc2}
\begin{split}
\Sigma_{(-1)}(R_1\Phi_n,1)
&=	q^{-n}\Sigma_{(-1)}(R_1\Phi,1),\\ 
\Sigma_{(0)}(R_1\Phi_n,1)
&	=q^{-n}\Sigma_{(0)}(R_1\Phi,1)-nq^{-n}\Sigma_{(-1)}(R_1\Phi,1),\\
\Sigma_{(-1)}(R_1\widehat{\Phi_n},1)
&	=q^{-2n}\Sigma_{(-1)}(R_1\hat\Phi,1),\\
\Sigma_{(0)}(R_1\widehat{\Phi_n},1)
&	=q^{-2n}\Sigma_{(0)}(R_1\hat\Phi,1)+nq^{-2n}\Sigma_{(-1)}(R_1\hat\Phi,1).
\end{split}
\end{align}
Hence, By \eqref{eq:calc1} and \eqref{eq:calc2},
we obtain the first three terms of
\eqref{eq:II}.
To get the last term of \eqref{eq:II},
note that simple calculation shows
\begin{align}
\begin{split}
\frac{1}{2}T^0(R_2\Phi,w)+\sum_{n\geq1}q^{-ns}T^0(R_2\Phi_{-n},w)
	&=T^-(R_2\Phi,s,w),\\
\frac{1}{2}T^0(R_2\hat\Phi,w)+\sum_{n\geq1}q^{-ns}T^0(R_2\widehat{\Phi_{-n}},w)
	&=T^+(R_2\hat\Phi,3-s,w).
\end{split}
\end{align}

On the other hand, by Lemma \ref{lem:T+rat}, \ref{lem:Trat},
$\tilde T^+(R_2\hat\Phi,3-s)$ and
$\tilde T^-(R_2\Phi,s)$
are rational functions of $q^s$,
and hence $Z(\Phi,s)$ is a rational function of $q^s$.
This finishes the proof.
\end{proof}

We define
\begin{equation}
\esR_1
=\frac{1}{\log q}\frac{2}{\rho_0}
=\frac{1}{\log q}\frac{2\zeta_k(2)}{q^{2-2\genus}\gC_k}.
\end{equation}
Together with Theorem \ref{th:gz} and Lemma \ref{lem:pd},
we obtain the following.
	\begin{cor}\label{cor:gz}
$Z(\Phi,s)$ is a rational function of $q^s$,
and $(1-q^{3-s})Z(\Phi,s)$ is holomorphic in the region $\res>2$.
Moreover,
\begin{equation*}
\Res_{s=3}Z(\Phi,s)=\esR_1\hat\Phi(0).
\end{equation*}
\end{cor}

We have the following functional equation for
the global zeta function $Z(\Phi,s)$.
\begin{cor}
Let
	$Z_{ad}(\Phi,s)
	=Z(\Phi,s)-\tilde T(R_2\Phi,s)$,
then $Z_{ad}(\Phi,s)$ satisfies the functional equation
\begin{equation*}
Z_{ad}(\Phi,s)=Z_{ad}(\hat\Phi,3-s).
\end{equation*}
\end{cor}

\section{Local theory}\label{sec:l}

\subsection{The canonical measure on the stabilizer}
\label{subsec:msr}
For an algebraic group $G$, Let $G^\circ$ denote its identity component.
In this subsection, we normalize invariant measures on 
$H_{k_v},H_{k_v}/{H_x^\circ}_{k_v},{H_x^\circ}_{k_v}$ for $x\in V_{k_v}'$,
and $H_\adl/{H_x^\circ}_\A, {H_x^\circ}_\adl$ for $x\in V_k''$,
following the method of \cite{kayu}.
We also define a constant $b_{x,v}$,
and compute the volume of ${H_x^\circ}_{\A}/{H_x^\circ}_{k}$
with respect to this measure (Proposition \ref{prop:mux}).

We define the invariant measure $dh_v$ on $H_{k_v}$
similarly as in \S\ref{subsec:dh1}.
Let $\K_{v}=G_{\co_v}$
be the standard maximal compact subgroup of $G_{k_v}$.
For $h_v\in H_{k_v}$, let
$h_v=\varrho(\kappa_v d(t_v,1)a(\tau_v)n(u_v))$
be its Iwasawa decomposition.
Define an invariant measure $dh_v$ on $H_{k_v}$ by
$dh_v=|\tau|_vd\kappa_v d^\times t_v d^\times\tau_v du_v$.
This normalization is equivalent to
$\int_{\varrho(\K_v)}dh_v=1$.

If we write $d_{\text{pr}}h=\vprod dh_v$, then by \eqref{productmeasures}
we can see
\begin{equation}\label{eq:dh/dprh}
dh=q^{1-\genus}\gC_k^{-2}d_{\text{pr}}h,
\end{equation}
where $dh$ is defined in \S\ref{subsec:dh1}.

Next, we choose a left invariant measure $dh_{x,v}'$
on $H_{k_v}/{H_x^\circ}_{k_v}$ for $x\in V_{k_v}'$.
Let $dy_v$ be the Haar measure on $V_{k_v}$ such that
the volume of $V_{\co_v}$ is one.
Note that $\frac{dy_v}{|P(y_v)|_v^{3/2}}$ is a
left $H_{k_v}$-invariant measure on $V_{k_v}'$ and
$H_{k_v}/{H_x^\circ}_{k_v}$ is a double cover of $H_{k_v}x\subset V_{k_v}'$.
Therefore there exists a unique left $H_{k_v}$-invariant measure
$dh_{x,v}'$ on $H_{k_v}/{H_x^\circ}_{k_v}$ such that
for $\Psi\in L^1(H_{k_v}x,\frac{dy_v}{|P(y_v)|_v^{3/2}})$,
\begin{equation}\label{eq:dh'}
\int_{H_{k_v}x}\Psi(y_v)\frac{dy_v}{|P(y_v)|_v^{3/2}}=
\int_{H_{k_v}/{H_x^\circ}_{k_v}}\Psi(h_{x,v}'\cdot x)dh_{x,v}'.
\end{equation}

If $h_{xy}\in H_{k_v}$ satisfies
$y=h_{xy}x$ and $i_{h_{xy}}$ is the inner automorphism 
$h\mapsto h_{xy}^{-1}hh_{xy}$ of $H_{k_v}$, then
$i_{h_{xy}}(H_{y\,k_v}^{\circ})=H_{x\,k_v}^{\circ}$ and hence $i_{h_{xy}}$
induces the map $i_{h_{xy}}: H_{k_v}/ H_{y\,k_v}^{\circ}\to
H_{k_v}/H_{x\,k_v}^{\circ}$. Since the integral on the right hand
side of (\ref{eq:dh'}) depends only on the orbit of $x$, it follows
that $i_{h_{xy}}^*(dh_{x,v}')=dh_{y,v}'$.

We divide $V_{k_v}'$ into three subsets
for conveniences.
Let
$$V_{k_v}'=V_{k_v}^{\rm sp}\amalg V_{k_v}^{\rm ur}\amalg V_{k_v}^{\rm rm},$$
where each subset consists of orbits corresponding to $k_v$,
quadratic unramified extension of $k_v$,
and quadratic ramified extension of $k_v$, respectively.
Only $V_{k_v}^{\rm rm}$ has two orbits and
$V_{k_v}^{\rm sp},V_{k_v}^{\rm ur}$ has a single orbit.

If we define
$$K_{x,k_v}=
\begin{cases}
k_v^\times\times k_v^\times	& x\in V_{k_v}^{\rm sp},\\
k_v(x)^\times 				& x\in V_{k_v}^{\rm ur},V_{k_v}^{\rm rm},\\
\end{cases}
$$
then by Proposition \ref{prop:Gx0},
${G_x^\circ}_{k_v}\cong K_{x,k_v}$.
Let $\co_{k_v(x)}$ denote the ring of integers of $k_v(x)$.
We will normalize the measure $dk_{x,v}$ on $K_{x,k_v}$ so that
$$\int_{\co_v^\times\times\co_v^\times}dk_{x,v}=1,\qquad
\int_{\co_{k_v(x)}^\times}dk_{x,v}=1$$
for each case.
This induces a Haar measure $dg_{x,v}''$ on ${G_x^\circ}_{k_v}$
via the isomorphism ${G_x^\circ}_{k_v}\cong K_{x,k_v}$.
Though there are two isomorphisms 
$\psi_{x,v}^{(1)},\psi_{x,v}^{(2)}$
for each case (see the proof of Proposition \ref{prop:Gx0}),
since
\begin{align*}
\psi_{x,v}^{(1)-1}(\co_v^\times\times\co_v^\times)
&	=\psi_{x,v}^{(2)-1}(\co_v^\times\times\co_v^\times)
&& x\in V_{k_v}^{\rm sp},\\
\psi_{x,v}^{(1)-1}(\co_{k_v(x)}^\times)
&	=\psi_{x,v}^{(2)-1}(\co_{k_v(x)}^\times)
&& x\in V_{k_v}^{\rm ur},V_{k_v}^{\rm rm},
\end{align*}
we can define $dg_{x,v}''$ without ambiguity.
Then, define $dh_{x,v}''$ by setting $dh_{x,v}''d^\times t_{1v}=dg_{x,v}''$
via the isomorphism ${H_x^\circ}_{k_v}\cong {G_x^\circ}_{k_v}/T_{k_v}$.

The next proposition shows that $dh_{x,v}''$
satisfies the functorial property.
In this sense, our choice $dh_{x,v}''$ is canonical.
\begin{prop}\label{prop:canonicity}
Suppose that $x,y\in V_{k_v}'$ and that $y=\varrho(g_{xy})x$
for some $g_{xy}\in G_{k_v}$. 
Let $i_{g_{xy}}:G_{y\,k_v}^{\circ}\to G_{x\,k_v}^{\circ}$
be the isomorphism given by $i_{g_{xy}}(g)=g_{xy}^{-1}gg_{xy}$.
Then,
\begin{equation}
dg_{y,v}''=i_{g_{xy}}^*(dg_{x,v}'')
\qquad \text{and}\qquad
dh_{y,v}''=i_{\varrho(g_{xy})}^*(dh_{x,v}'')\,.
\end{equation}
\end{prop}
\proof
We only consider the case $[k_v(x):k_v]=2$.
We can prove the case $k_v(x)=k_v$ similarly.
One can easily show that for any $g_{xy}$,
there exist $\psi_{x,v},\psi_{y,v}$ such that
the following diagram is commutative.
$$
\begin{CD}
{G_y^\circ}_{k_v}	@>{\psi_{y,v}}>>	k_v(y)^\times	\\
@V{g_{xy}}VV							@|				\\
{G_x^\circ}_{k_v}	@>{\psi_{x,v}}>>	k_v(x)^\times	\\
\end{CD}
$$
This establishes the first claim and the second claim follows from the
observation that $i_{g_{xy}}|_{T_{k_v}}$ is the identity map.
\hspace*{\fill}$\square$

Define a constant $b_{x,v}>0$ such that
$dh_v=b_{x,v}dh_{x,v}'dh_{x,v}''.$
Then, the following proposition shows that
$b_{x,v}$ depends only on the orbit of $x$.

\begin{prop}\label{bxindep} 
If $x,y\in V'_{k_v}$ are in the same $H_{k_v}$-orbit,
then $b_{x,v}=b_{y,v}$.  
\end{prop} 
\proof
 Since the group $H_{k_v}$ is unimodular,
$i_{h_{x,y}}^* dh_v = dh_v$.  
Hence,
\begin{equation*}
\begin{aligned}
dh_v & = b_{y,v} dh_{y,v}' dh''_{y,v} \\
& = b_{y,v} i_{h_{x,y}}^* dh_{x,v}'i_{h_{x,y}}^* dh_{x,v}''
& = b_{y,v}b_{x,v}^{-1} i_{h_{x,y}}^* dh_v \\
& = b_{y,v}b_{x,v}^{-1} dh_v\,. 
\end{aligned}
\end{equation*}
Therefore $b_{x,v}=b_{y,v}$.
\hspace*{\fill}$\square$

For $x\in V_k''$,
Let
$$
dh_x'=\vprod b_{x,v}dh_{x,v}',\qquad
dh_x''=\vprod dh_{x,v}''
$$
be measures on $H_\adl/{H_x^\circ}_\A, {H_x^\circ}_\adl$, respectively.
Then,
$$dh_x'dh_x''=d_{\text{pr}}h.$$

We will conclude this subsection by computing
the volume of ${H_x^\circ}_{\A}/{H_x^\circ}_{k}$
with respect to the measure $dh_x''$ defined above.

\begin{prop}\label{prop:mux}
For $x\in V_k''$,
$$\int_{(H_x^\circ)_\A/(H_x^\circ)_k}dh_x''
=2\frac{\gC_{k(x)}}{\gC_k}.$$
\end{prop}

\proof
Recall that 
$(H_x^\circ)_\A/(H_x^\circ)_k\cong(G_x^\circ)_\A/T_\A(G_x^\circ)_k$.
One can easily see that the inclusion
$$
\xymatrix@1
{(G_x^\circ)^0_\A/T_\A^0(G_x^\circ)_k	\ar@{^{(}->}[r]&
(G_x^\circ)_\A/T_\A(G_x^\circ)_k}
$$
has index two,
and the sequence
$$
\begin{CD}
1				@>>>
T_\A^0/T_k		@>>>
(G_x^\circ)^0_\A/(G_x^\circ)_k		@>>>
(G_x^\circ)^0_\A/T_\A^0(G_x^\circ)_k	@>>>
1
\end{CD}
$$
is exact.
Since
$$
\int_{T_\A^0/T_k}d^\times t=\gC_k\quad{\rm and}\quad
\int_{(G_x^\circ)^0_\A/(G_x^\circ)_k}dg_x''=\gC_{k(x)},
$$
we can get
$$
\int_{(H_x^\circ)_\A/(H_x^\circ)_k}dh_x''
=2\int_{(G_x^\circ)^0_\A/T_\A^0(G_x^\circ)_k} dh_x''
=2\frac{\gC_{k(x)}}{\gC_k}.
$$
\hspace*{\fill}$\square$

\subsection{Local zeta function}
Here, we will define the local zeta function,
and compute the local zeta function
and the constant $b_{x,v}$ for some orbital representatives.
\begin{defn}
Let $x\in V_{k_v}'$,
$\Phi_v\in\cs(V_{k_v})$,
and $s\in\C$.
We define
\begin{align*}
Z_{x,v}(\Phi_v,s)
&	=b_{x,v}\int_{H_{k_v}/{H_x^\circ}_{k_v}}
		|\chi(h_{x,v}')|_v^s\Phi_v(h_{x,v}'x)dh_{x,v}',\\
\Omega_{x,v}(\Phi_v,s)
&	=\int_{H_{k_v}x}|P(y)|_v^{s/2}\Phi_v(y)\frac{dy}{|P(y)|_v^{3/2}}.
\end{align*}
\end{defn}
The function
$\Omega_{x,v}(\Phi_v,s)$
is called the local zeta function. 
By the definition of $dh_{x,v}'$,
\begin{align*}
Z_{x,v}(\Phi_v,s)
&	=b_{x,v}\int_{H_{k_v}/{H_x^\circ}_{k_v}}
	\frac{|P(h_{x,v}'x)|_v^{s/2}}{|P(x)|_v^{s/2}}
		\Phi_v(h_{x,v}'x)dh_{x,v}'\\
&	=b_{x,v}|P(x)|_v^{-s/2}\Omega_{x,v}(\Phi_v,s).
\end{align*}
Since $\Omega_{x,v}(\Phi_v,s)$ depends only on the orbit of $x$,
for $x,y\in V_{k_v}'$ in the same orbit,
\begin{equation}\label{eq:sameorbit}
Z_{x,v}(\Phi_v,s)=
\frac{|P(y)|_v^{s/2}}{|P(x)|_v^{s/2}}
Z_{y,v}(\Phi_v,s).
\end{equation}

Now, we will express  $Z_{x,v}$ and $b_{x,v}$ explicitly
for some representative $x$.
\begin{defn}
We call $w_v\in V_{k_v}'$ a standard orbital representative if
\begin{equation*}
F_{w_v}(z_1,z_2)=
	\begin{cases}
		z_1z_2 & {w_v}\in V_{k_v}^{\rm sp},\\
		(z_1+\theta z_2)(z_1+\theta' z_2)
			& {w_v}\in V_{k_v}^{\rm ur},V_{k_v}^{\rm rm},
			\text{where } \co_{k_v({w_v})}=\co_v[\theta].
	\end{cases}
\end{equation*}
For each orbit in $V_{k_v}'$,
we take one of the standard orbital representatives
and denote the fixed set of representatives by $\sr_v$.
\end{defn}
Note that for a standard orbital representative ${w_v}$,
$P({w_v})$ is the discriminant of $k_v(x)$ over $k_v$.
Hence,
\begin{equation}\label{eq:pw}
|P({w_v})|_v=
	\begin{cases}
		1 &	{w_v}\in V_{k_v}^{\rm sp},V_{k_v}^{\rm ur},\\
		q_v^{-1} & {w_v}\in V_{k_v}^{\rm rm}.
	\end{cases}
\end{equation}

Let $\Phi_{v,0}$ be the characteristic function of $V_{\co_v}$.
Firstly, we will give the explicit formula of
$Z_{{w_v},v}(\Phi_{v,0},s)$.
Although our choice of the measure on the stabilizers is
different from that of \cite{Dat1} in general,
they coincide for standard orbital representatives.
Hence we can use his result directly.
Note that our local zeta functions are $b_{x,v}$ times
that of Datskovsky's.

\begin{prop}[\cite{Dat1}Proposition 4.1]\label{prop:zx}
For a standard orbital representative ${w_v}$,
\begin{equation*}
Z_{{w_v},v}(\Phi_{v,0},s)=
	\begin{cases}
		\displaystyle{\frac{1}{1-q^{1-s}}} & {w_v}\in V_{k_v}^{\rm sp},\\
		\displaystyle{\frac{1+q^{-s}}{(1-q^{-s})(1-q^{1-s})}} 
				& {w_v}\in V_{k_v}^{\rm ur},\\
		\displaystyle{\frac{1}{(1-q^{-s})(1-q^{1-s})}}
				& {w_v}\in V_{k_v}^{\rm rm}.\\
	\end{cases}
\end{equation*}
\end{prop}

Secondly, we will give the values of $b_{{w_v},v}$.
To use Datskovsky's result, we need some discussion.
\begin{lem}
For a standard orbital representative ${w_v}$,
$$\int_{\varrho(\K_v)\cap {H_{w_v}^\circ}_{k_v}}dh_{{w_v},v}''=1.$$
\end{lem}
\proof
Here, we only consider the case $[k_v({w_v}):k_v]=2$.
The case $k_v({w_v})=k_v$ can be proved similarly.
Let $N=N_{k_v(w_v)/k_v}:k_v(w_v)^\times\to k_v^\times$ be the norm map.

Let $F_{w_v}(z_1,z_2)=(z_1+\theta z_2)(z_1+\theta'z_2)$.
Recall that for $g=\left(t,\stwtw abcd\right)\in G_{k_v}$,
the condition 
$g\in {G_{w_v}^\circ}_{k_v}$
is equivalent to
\begin{equation*}
tN(a+b\theta)=1,\quad c+d\theta=\theta(a+b\theta)
\end{equation*}
and that
the isomorphism
$\psi_{w_v,v}: {G_{w_v}^\circ}_{k_v}\rightarrow k_v(w_v)^\times$
is defined by
$g\mapsto a+b\theta$
for this $g$.
This map also gives the isomorphism
$T_{k_v}\cong k_v^\times$.
One can easily see that
$\psi_{w_v,v}(T_{k_v}\cap\ck_v)=\co_v^\times$.
We claim that 
$$
\psi_{w_v,v}({G_{w_v}^\circ}_{k_v}\cap\ck_v)={\co_{k_v(w_v)}}^\times.
$$
Let $g=(t,\stwtw abcd)\in{G_{w_v}^\circ}_{k_v}\cap\ck_v$.
Then $N(a+b\theta)=t^{-1}\in\co_v^\times$
and hence we have $a+b\theta\in {\co_{k_v(w_v)}}^\times$.
On the other hand,
any element of ${\co_{k_v(w_v)}}^\times$
can be written as
$\gamma=a+b\theta$ with $a, b\in\co_v$.
For this $\gamma$, take $t,c,d\in k_v$ such that
\begin{equation*}
tN(a+b\theta)=1,\quad c+d\theta=\theta(a+b\theta),
\end{equation*}
and let $g=(t,\stwtw abcd)\in G_{k_v}$.
Then, clearly
$t\in\co_v^\times, c,d\in\co_v$,
$g\in {G_{w_v}^\circ}_{k_v}$,
 and 
$\gamma=\psi_{w_v,v}(g)$.
Also,
$P(w_v)=P(g\cdot w_v)=\chi(g)^3P(w_v)$
shows $\chi(g)=t(ad-bc)\in\co_v^\times$.
Hence, we have $g\in\ck_v$.
This establishes the claim and now
by the definition of $dh_{{w_v},v}''$,
we have
$$
\int_{\varrho(\K_v\cap {G_{w_v}^\circ}_{k_v})}dh_{{w_v},v}''
=	\frac{{\rm vol}({\co_{k_v(w_v)}}^\times)}{{\rm vol}(\co_v^\times)}
=	\frac{1}{1}=1.
$$
\hspace*{\fill}$\square$

\begin{prop}\label{prop:bw}
For a standard orbital representative ${w_v}$,
$$b_{{w_v},v}=\frac{|P({w_v})|_v^{3/2}}{\int_{\varrho(\K_v){w_v}}dx_v}.$$
\end{prop}

\proof
Recall that $dh_v$ is the measure on $H_{k_v}$ such that
${\rm vol}(\varrho(\K_v))=1$. Hence,
\begin{align*}
1
=	\int_{\varrho(\K_v)}dh_v
&=	b_{{w_v},v}\int_{\varrho(\K_v)H_{w_v}^\circ/H_{w_v}^\circ}dh_{{w_v},v}'
	\int_{\varrho(G_{w_v}^\circ\cap\K_v)}dh_{{w_v},v}''\\
&=	b_{{w_v},v}\int_{\varrho(\K_v)H_{w_v}^\circ/H_{w_v}^\circ}dh_{{w_v},v}'\\
&=	b_{{w_v},v}\int_{\varrho(\K_v){w_v}}\frac{dx_v}{|P(x_v)|_v^{3/2}}.
\end{align*}
On the other hand, for $x_v\in\varrho(\K_v){w_v}, |P(x_v)|_v=|P({w_v})|_v$.
This finishes the proof.
\hspace*{\fill}$\square$

Under the above preparation, we can describe $b_{{w_v},v}$.
For our purpose, $|P({w_v})|_v^{3/2}/b_{{w_v},v}$
(which is equal to $\int_{\varrho(\K_v){w_v}}dx_v$
 by Proposition \ref{prop:bw})
is more important than $b_{{w_v},v}$ itself,
and the following proposition gives the desired description.
\begin{prop}\label{prop:bx}
For a standard orbital representative ${w_v}$,
\begin{equation*}
\frac{|P({w_v})|_v^{3/2}}{b_{{w_v},v}}=
	\begin{cases}
		\frac{1}{2}(1-q_v^{-2}) & {w_v}\in V_{k_v}^{\rm sp},\\
		\frac{1}{2}(1-q_v^{-1})^2 
				& {w_v}\in V_{k_v}^{\rm ur},\\
		\frac{1}{2}q_v^{-1}(1-q_v^{-1})(1-q_v^{-2})
				& {w_v}\in V_{k_v}^{\rm rm}.\\
	\end{cases}
\end{equation*}
\end{prop}
%
This proposition is essentially
proved in \cite{Dat1} Proposition 4.2,
but he did not give a
detailed account of the argument
for determining the volume of each of $\varrho(\ck_v)w_v$
from the sum for two orbits in $V_{k_v}^{\rm rm}$.
We can directly determine each volume
by following the argument in \cite{kayu2} Section 4.
For the convenience of the reader,
we briefly sketch their argument
by indicating the main steps of the proof.
Note that the action of
$\ck_v=G_{\co_v}$ on $V_{\co_v}$
induces the action of
$G_{\co_v/\gp_v^2}$ on $V_{\co_v/\gp_v^2}$.
Let $x$ be one of the standard orbital
representatives in $V_{k_v}^{\rm rm}$
and denote by
$\overline x\in V_{\co_v/\gp_v^2}$
its reduction modulo $\gp_v^2$.

\begin{enumerate}
\item
	Set
	$\cd=\{y\in V_{\co_v}\mid y\equiv x (\gp_v^2)\}$.
\item
	If $y\in\cd$, $k_v(y)=k_v(x)$.
	Moreover, $\cd\subset \ck_v x$.
\item
	$\vol(\ck_vx)=\vol(\cd)
		\#(G_{\co_v/\gp_v^2}/G_{\co_v/\gp_v^2,\overline x})$.
\item
	$\vol(\cd)=q_v^{-6}$,
	$\#G_{\co_v/\gp_v^2}=q_v^6(q_v-1)^2(q_v^2-1)$,
\item
	$\#G_{\co_v/\gp_v^2,\overline x}=2(q_v-1)q_v^4$.
\end{enumerate}
From this proposition, we can obtain $b_{{w_v},v}$ easily using
\eqref{eq:pw}.

\section{The mean value theorem}\label{sec:mv}
In this section, we will deduce our mean value theorem
by putting together the results we have obtained before.
In \S\ref{subsec:as}, we will see that
the global zeta function is approximately
the Dirichlet generating series
for the sequence ${\rm vol}(H_{x\A}^\circ/H_{xk}^\circ)$.
If it were exactly this generating series,
the theory of partial fraction would allow us
to extract the mean value of the coefficients
from the analytic behavior of this series.
However, our global zeta function contains
an additional factor in each term.
In \S\ref{subsec:fp} we will surmount this difficulty
by using the technique called the filtering process,
which was formulated by Datskovsky-Wright \cite{Da-Wr1}.

\subsection{The adelic synthesis}\label{subsec:as}
We will introduce some notation.
For the rest of this paper, we suppose $\Phi\in\cs(V_\A)$ is of the form
$\Phi=\otimes\Phi_v$, where $\Phi_v\in\cs(V_{k_v})$.
For $x\in V_k''$, define $Z_x(\Phi,s)=\vprod Z_{x,v}(\Phi_v,s)$.
For each $v\in\gM$,
take $w_{v,x}\in\sr_v$
which lies in the orbit of $x$.
We write
$\Xi_{x,v}(\Phi_v,s)= Z_{{w_{v,x}},v}(\Phi_v,s)$ and 
$\Xi_{x}(\Phi,s)=\vprod\Xi_{x,v}(\Phi_v,s)$. 

Then, as is well known,
our global zeta function has the following expansion.
(See \cite{kayu} Section 6, for example.
Note that $|H_x/H_x^\circ|=2$ for all $x\in V_k''$.)
\begin{align*}\label{eq:azd}
Z(\Phi,s)
&=	\frac{g^{1-\genus}}{2\gC_k^2}\sum_{x\in H_k\backslash V_k''}
	\int_{(H_x^\circ)_\A/(H_x^\circ)_k}dh_x''
	\int_{H_\adl/(H_x^\circ)_\adl}
	\idn{\chi(h_x')}^s\Phi(h_x'\cdot x)dh_x'\\
&=	q^{1-\genus}\gC_k^{-3}
		\sum_{x\in H_k\backslash V_k''}\gC_{k(x)}Z_x(\Phi,s).
\end{align*}
We will consider $Z_x(\Phi,s)$.
By \eqref{eq:sameorbit}, we have
\begin{align*}
Z_x(\Phi,s)
&	=\vprod Z_{x,v}(\Phi_v,s)
	=\vprod \frac{|P(w_{v,x})|_v^{s/2}}{|P(x)|_v^{s/2}}
				\Xi_{x,v}(\Phi_v,s).
\end{align*}
Observe that since $P(x)\in k^\times$,
$\vprod|P(x)|_v=\idn{P(x)}=1$.
Also since $P(w_v)$ is the local discriminant of $k_v(x)$ over $k_v$,
$\vprod|P(w_{v,x})|_v=\nd{k(x)}^{-1}$.
Hence, we obtain the following:
\begin{equation}
Z(\Phi,s)
	=q^{1-\genus}\gC_k^{-3}
		\sum_{x\in H_k\backslash V_k''}
			\frac{\gC_{k(x)}}{\nd{k(x)}^{s/2}}\Xi_x(\Phi,s).
\end{equation}

Let $T$ be a finite set of places of $k$,
and we denote by $\sr_T$ the Cartesian product $\prod_{v\in T}\sr_v$.
We consider $T$-tuples $w_T=(w_v)_{v\in T}$
of an element of $\sr_T$.
We say that $x\in V_k''$ is equivalent to $w_T$
if $x$ lies in the $H_{k_v}$-orbit of $w_v$ for each $v\in T$
and denoted by $x\sim w_T$.
Let $\Phi|_T=\prod_{v\in T}\Phi_v$, and
\begin{equation}\label{eq:zxs}
Z_{w_T,T}(\Phi|_T,s)
	=\prod_{v\in T}Z_{w_v,v}(\Phi_v,s).
\end{equation}
Then,
\begin{equation}\label{eq:ccccc}
Z(\Phi,s)
	=q^{1-\genus}\gC_k^{-3}
		\sum_{w_T\in\sr_T}Z_{w_T,T}(\Phi|_T,s)
\left(
		\sum_{x\sim w_T}\frac{\gC_{k(x)}}{\nd{k(x)}^{s/2}}
		\prod_{v\not\in T}\Xi_{x,v}(\Phi_v,s)
\right).
\end{equation}
For $v\not\in T$, we take
$\Phi_v$ as the characteristic function $\Phi_{v,0}$
of $V_{\co_v}$.
Then, for $v\not\in T$, $\Xi_{x,v}(\Phi_{v,0},s)$ is given in
Proposition \ref{prop:zx} and hence we have
\begin{equation}\label{eq:thl}
\prod_{v\not\in T}\Xi_{x,v}(\Phi_{v,0},s)
	=\frac{\zeta_{k,T}(s-1)\zeta^2_{k,T}(s)}
		{\zeta_{k(x),T}(s)},
\end{equation}
where $\zeta_{k,T}(s)$ and $\zeta_{k(x),T}(s)$ are
the truncated Dedekind zeta function
\begin{equation*}
\zeta_{k,T}(s)=\prod_{v\not\in T}(1-q_v^{-s})^{-1},\quad
\zeta_{k(x),T}(s)=\prod_{\substack{\mu\in\gM_{k(x)}\\ \mu|v,v\not\in T}}
					(1-q_\mu^{-s})^{-1}.
\end{equation*}
We will denote the function \eqref{eq:thl} by $\eta_{k(x),T}(s)$.
Note that $\eta_{k(x),T}(s)$ is
a Dirichlet series.
Set the Dirichlet series $\xi_{w_T}(s)$ by
\begin{equation}
\xi_{w_T}(s)=
	\sum_{\substack{[k(x):k]=2\\ x\sim w_T}}
		\frac{\gC_{k(x)}}{\nd{k(x)}^{s/2}}
		\eta_{k(x),T}(s).
\end{equation}
Then \eqref{eq:ccccc} becomes
\begin{equation}\label{eq:dsz}
Z(\Phi,s)
	=q^{1-\genus}\gC_k^{-3}
		\sum_{w_T\in\sr_T}Z_{w_T}(\Phi|_T,s)\xi_{w_T}(s).
\end{equation}
In order to determine the analytic properties
of $\xi_{w_T}(s)$, we require the following lemma.
This is quite similar to \cite{kayu} Lemma 6.17,
and we omit the proof.
\begin{lem}\label{lem:choosephi}
Let $v\in\gM, x\in V_{k_v}'$ and $s_0\in\C$.
Then there exists $\Phi_v\in\cs(V_{k_v})$
such that the support of $\Phi_v$
is contained in $H_{k_v}x$,
$Z_{x,v}(\Phi_v,s)$ is a polynomial in $q_v^s, q_v^{-s}$
and $Z_{x,v}(\Phi_v,s_0)\not=0$.
\end{lem}
Let
\begin{equation*}
\esR_2=\frac{2\zeta_k(2)\gC_k^2}{\log q}.
\end{equation*}
Also for $w_v\in\sr_v$ and $w_T=(w_v)_{v\in T}\in\sr_T$, we define 
\begin{equation*}
\vep_v(w_v)=\frac{|P(w_v)|_v^{3/2}}{b_{w_v}},\qquad
\vep_T(w_T)=\prod_{v\in T}\vep_v(w_v).
\end{equation*}
Now we can prove the following theorem.
\begin{thm}\label{th:xi}
For $w_T=(w_v)_{v\in T}$,
	the series $\xi_{w_T}(s)$ is a rational function of $q^s$,
	and is holomorphic in the region $\res>3$.
	Also $(1-q^{3-s})\xi_{w_T}(s)$ is holomorphic in the region
	$\res>2$.
Moreover,
	\begin{equation*}
		\Res_{s=3}\xi_{w_T}(s)
		=\esR_2\vep_T(w_T)
	\end{equation*}
\end{thm}
\proof \ 
For each $v\in T$, we choose $\Phi_v\in\cs(V_{k_v})$
such that $\supp(\Phi_v)\subset H_{k_v}w_v$.
Then, $Z_{w_v,v}(\Phi_v,s)=0$ unless $x\sim w_v$.
Hence \eqref{eq:dsz} becomes
\begin{equation}
Z(\Phi,s)
	=q^{1-\genus}\gC_k^{-3}
		Z_{w_T}(\Phi|_T,s)\xi_{w_T}(s).
\end{equation}
Then the first two statements follows from Corollary \ref{cor:gz}
 and Lemma \ref{lem:choosephi}.
We will compute $\Res_{s=3}\xi_{w_T}(s)$.
By Corollary \ref{cor:gz},
$\Res_{s=3}Z(\Phi,s)=
\esR_1\hat\Phi(0)$.
We consider
\begin{equation*}
\hat\Phi(0)=\int_{V_\adl}\Phi(x)dx=q^{3-3\genus}
\vprod\int_{V_{k_v}}\Phi_v(x_v)dx_v.
\end{equation*}
For $v\not\in T$, 
$\int_{V_{k_v}}\Phi_v(x_v)dx_v=1$ since
$\Phi_v=\Phi_{v,0}$ is the characteristic function of $V_{\co_v}$.
For $v\in T$,
\begin{equation*}
\int_{V_{k_v}}\Phi_v(x_v)dx_v
=	\int_{H_{k_v}w_v}\Phi_v(x_v)dx_v
=	\Omega_{w_v,v}(\Phi_v,3)
=	\frac{|P(w_v)|_v^{3/2}}{b_{w_v}}
		Z_{w_v,v}(\Phi_v,3).
\end{equation*}
Hence we have
\begin{equation*}
\hat\Phi(0)=q^{3-3\genus}
	Z_{w_T}(\Phi|_T,3)\vep_T(w_T).
\end{equation*}
Together with \eqref{eq:dsz}, this yields the residue of 
$\xi_{w_T}(s)$.
\hspace*{\fill}$\square$

\subsection{The filtering process}\label{subsec:fp}

We fix a finite set $T_0$ of places of $k$ and
$w_{T_0}=(w_v)_{v\in{T_0}}$.
\begin{defn}
For each finite subset $T\supset T_0$ of $\gM$,
we define
\begin{equation*}
\xi_{w_{T_0},T}(s)
=	\sum_{x\sim w_{T_0}}
		\frac{\gC_{k(x)}}{\nd{k(x)}^{s/2}}\eta_{k(x),T}(s).
\end{equation*}
\end{defn}
For $v\in\gM$, let
$$
E_v
=\sum_{w_v\in\sr_v}\vep_v(w_v)
=	1-q_v^{-2}-q_v^{-3}+q_v^{-4},
$$
and also for any subset $T'$ of $\gM$,
define
$$
E_{T'}=\prod_{v\in T'}E_v.
$$
Note that this product always converges to a positive number.
\begin{prop}
The Dirichlet series $\xi_{w_{T_0},T}(s)$
becomes a rational function of $q^s$ and holomorphic in the region $\res>3$.
Also $(1-q^{3-s})\xi_{w_{T_0},T}(s)$
is holomorphic in the region $\res>2$.
The residue of $\xi_{w_{T_0},T}(s)$ at $s=3$ is given by
\begin{equation*}
	\esR_2\vep_{T_0}(w_{T_0})E_{T\backslash T_0}.
\end{equation*}
\end{prop}
\begin{proof}
For $y_T=(y_v)_{v\in T}\in\sr_T$,
we denote $y_T|_{T_0}=(y_v)_{v\in T_0}\in\sr_{T_0}$.
Then
\begin{equation*}
\xi_{w_{T_0},T}(s)=\sum_{y_T|_{T_0}=w_{T_0}}\xi_{y_T}(s).
\end{equation*}
Now the proposition immediately follows from Theorem \ref{th:xi}.
\end{proof}

To deduce our mean value theorem,
we have to show some properties of $\eta_{k(x),T}(s)$.
For two Dirichlet series
$\vartheta_i(s)=\sum_{n=0}^{\infty}r_{i,n}/q^{ns}, i=1,2$,
we will indicate $\vartheta_1(s)\preceq \vartheta_2(s)$
or $\vartheta_2(s)\succeq \vartheta_1(s)$
if $r_{1,n}\leq r_{2,n}$ for all $n$.
Especially, write $\vartheta(s)=\sum_{n=0}^{\infty}r_n/q^{ns}\succeq0$
if $r_n\geq0$ for all $n$.
The following proposition is easy to prove.
\begin{prop}\label{prop:sd}
The Dirichlet series $\eta_{k(x),T}(s)$
satisfies $\eta_{k(x),T}(s)\succeq0$,
and its first coefficient is $1$.
Also, for all $k(x)$,
\begin{equation*}
\eta_{k(x),T}(s)\preceq\eta_{T}(s)
	=\frac{\zeta_{k,{T}}(s-1)\zeta_{k,{T}}^2(s)}{\zeta_{k,{T}}(2s)}.
\end{equation*}
Moreover, $\eta_{T}(s)$ converges in the region
$\res>2$ and 
\begin{equation*}
\lim_{T\uparrow \gM}(\eta_T(1)-1)=0.
\end{equation*}
\end{prop}
Let us define $a_n\geq0$ by
\begin{equation*}
\sum_{n\geq0}\frac{a_n}{q^{ns}}
=\sum_{x\sim w_{T_0}}\frac{\gC_k(x)}{\nd{k(x)}^{s/2}}.
\end{equation*}
Now, we are ready to prove the following theorem.
\begin{thm}\label{thm:pmv}
\begin{equation*}
\lim_{n\to\infty}\frac{a_n}{q^{3n}}
=\log q\esR_2\vep_{T_0}(w_{T_0})E_{\gM\setminus {T_0}}.
\end{equation*}
\end{thm}
\proof
Since $\eta_{k(x),T}(s)\succeq0$, we have
$$
\xi_{w_{T_0},T}(s)\succeq
\sum_{x\sim w_{T_0}}\frac{\gC_k(x)}{\nd{k(x)}^{s/2}}
=\sum_{n\geq0}\frac{a_n}{q^{ns}}.
$$
Hence, if one write $\xi_{w_{T_0},T}(s)=\sum{r_{T,n}}/{q^{ns}}$,
then $r_{T,n}\geq a_n$.
By the theory of partial fraction,
$\lim_{n\to\infty}{r_{T,n}}/{q^{3n}}
=\log q \esR_2\vep_{T_0}(w_{T_0})E_{T\setminus {T_0}}$.
Hence,
$$
\overline{\rm lim}_{n\to\infty}
\frac{a_n}{q^{3n}}\leq
\log q \esR_2\vep_{T_0}(w_{T_0})E_{T\setminus {T_0}}.
$$
By letting $T$ approach to $\gM$, we obtain
$$
\overline{\rm lim}_{n\to\infty}\frac{a_n}{q^{3n}}
\leq \log q\esR_2 \vep_{T_0}(w_{T_0})E_{\gM\setminus {T_0}}.
$$

This allows us to take $R'>0$ such that 
$a_n\leq q^{3n}R'$ for all $n$.
Let $\eta_{T}(s)=\sum_{n\geq0}l_{T,n}/q^{ns}$. Then $l_{T,0}=1$ and
$$
\xi_{w_{T_0},T}(s)
\preceq		\sum_{n\geq0}\frac{a_n}{q^{ns}}\eta_{T_0}(s)
=			\sum_{n\geq0}\frac{\sum_{n_1+n_2=n}
				a_{n_1}l_{T,n_2}}{q^{ns}}.
$$
Since
\begin{align*}
\sum_{n_1+n_2=n}a_{n_1}l_{T,n_2}
&=		a_n+\sum_{n_2=1}^nl_{T,n_2}a_{n-n_2}\\
&\leq	a_n+q^{3n}R'\sum_{n_2=1}^n\frac{l_{T,n_2}}{q^{3n_2}}
\leq	a_n+q^{3n}R'(\eta_{T}(1)-1),
\end{align*}
we have
$$
\underline{\rm lim}_{n\to\infty}\frac{a_n}{q^{3n}}
\geq \log q \esR_2\vep_{T_0}(w_{T_0})E_{T\setminus {T_0}}-R'(\eta_T(1)-1).
$$
Again by letting $T$ approach to $\gM$, we obtain
$$
\underline{\rm lim}_{n\to\infty}\frac{a_n}{q^{3n}}
\geq \log q \esR_2\vep_{T_0}(w_{T_0})E_{\gM\setminus {T_0}}.
$$
Together with the estimate for the superior limit,
we obtain the result.
\hspace*{\fill}$\square$

\subsection{Main results}\label{subsec:mr}
Let us rewrite Theorem \ref{thm:pmv} to a mean value theorem
for the degree zero divisor class groups of quadratic extensions.
Let $\sB_v$ be the index set
of extensions of $k_v$ of degree not greater than two.
By assumption that ${\rm char}(k)\not=2$, the cardinality of
this set is four for all $v$.
We denote by $k_{v,\beta_v}$ an extension corresponding to $\beta_v\in\sB_v$.
From now on, the letter $L$ always denotes a quadratic extension of $k$.
For $\beta_v\in\sB_v$,
we write $L\sim\beta_v$
when the extension of $L/k$ at $v$
is $k_{v,\beta_v}/k_v$.
We fix a finite set $T$ of places of $k$ and
$\beta_T=(\beta_v)_{v\in T}\in \prod_{v\in T}\sB_v$.
If $L\sim\beta_v$ for all $v\in T$
then we write $L\sim\beta_T$.

Define $b_{\beta_v}$ and $b_{\beta_T}$ as follows.
\begin{equation*}
\begin{split}
& b_{\beta_v}=
\begin{cases}
\tfrac12(1-q_v^{-2}) & k_{v,\beta_v}=k_v,\\
\tfrac12(1-q_v^{-1})^2 &
	k_{v,\beta_v}\text{ is quadratic unramified over }k_v, \\
\tfrac12q_v^{-1}(1-q_v^{-1})(1-q_v^{-2}) &
	k_{v,\beta_v}\text{ is quadratic ramified over }k_v, \\
\end{cases}\\
& b_{\beta_T}=\prod_{v\in T}b_{\beta_v}.\\
\end{split}
\end{equation*}

For all quadratic extensions $L$
except $k\otimes_{\F_q}\F_{q^2}$, $q_L=q$.
Then, we can rewrite Theorem \ref{thm:pmv} as follows.
\begin{thm}\label{thm:mv}
\begin{equation*}
\lim_{n\to\infty}\frac{1}{q^{3n}}
\sum_{\substack{L\sim\beta_S\\ \nd{L}=q^{2n}}}h_{L}
=	2\gC_k h_k\zeta_k(2)b_{\beta_S}
	\prod_{v\not\in S}(1-q_v^{-2}-q_v^{-3}+q_v^{-4}).
\end{equation*}
\end{thm}
If we take $S=\emptyset$, we will obtain
Theorem \ref{thm:intro} in the introduction.

We conclude this paper
with some modification of this formula.
The next proposition
about the density of quadratic extensions
is well known.
We can evaluate this formula
by means of class field theory
or an analysis of a slight variation
of Tate's zeta function.
We simply state the result here.

\begin{prop}
Set
\begin{equation*}
c_{\beta_S}=\prod_{v\in S}c_{\beta_v},\quad
c_{\beta_v}=
\begin{cases}
\tfrac12(1-q_v^{-1}) &
	k_{v,\beta_v}\text{ is unramified over }k_v, \\
\tfrac12q_v^{-1}(1-q_v^{-1}) &
	k_{v,\beta_v}\text{ is quadratic ramified over }k_v. \\
\end{cases}
\end{equation*}
Then
$$\lim_{n\to\infty}\frac{1}{q^{2n}}
\sum_{\substack{L\sim\beta_S\\ \nd{L}=q^{2n}}}1
=	2q^{1-\genus}\gC_k c_{\beta_S}
	\prod_{v\not\in S}(1-q_v^{-2}).$$
\end{prop}

Therefore, together with Theorem \ref{thm:mv},
we can obtain the following formula:
\begin{cor}
\begin{equation*}
\lim_{n\to\infty}\frac{1}{q^{n}}
\frac{\sum_{L\sim\beta_S,\nd{L}=q^{2n}}h_{L}}
{\sum_{L\sim\beta_S,\nd{L}=q^{2n}}1}
=	q^{\genus-1}h_k\zeta_k(2)d_{\beta_S}
	\prod_{v\not\in S}\frac{1+q_v^{-1}-q_v^{-3}}{1+q_v^{-1}},
\end{equation*}
where, $d_{\beta_S}$ is given by
\begin{equation*}
d_{\beta_S}=\prod_{v\in S}d_{\beta_v},\qquad
d_{\beta_v}=
\begin{cases}
1+q_v^{-1} & k_{v,\beta_v}=k_v,\\
1-q_v^{-1} &
	k_{v,\beta_v}\text{ is quadratic unramified over }k_v, \\
1-q_v^{-2} &
	k_{v,\beta_v}\text{ is quadratic ramified over }k_v. \\
\end{cases}
\end{equation*}
\end{cor}
Let $\genus_L$ be the genus of $L$.
To avoid the notational confusion,
here we denote the genus of $k$ by $\genus_k$.
If $\nd{L}=q^{2n}$, $\genus_L-1=2\genus_k-2+n$.
Hence the preceding formula can also be expressed as follows:
\begin{cor}
\begin{equation*}
\lim_{n\to\infty}\frac{1}{q^{n}}
\frac{\sum_{L\sim\beta_S,\genus_L=n}h_{L}}
{\sum_{L\sim\beta_S,\genus_L=n}1}
=	\frac{h_k\zeta_k(2)}{q^{\genus_k}}d_{\beta_S}
	\prod_{v\not\in S}\frac{1+q_v^{-1}-q_v^{-3}}{1+q_v^{-1}}.
\end{equation*}
\end{cor}

\end{document}